%% file: Encrustamientoscirculares.tex
\documentclass[11pt]{article}
\usepackage{latexsym,xy,eucal,mathrsfs,graphs}
\usepackage{setspace}
\usepackage{amsthm}
\usepackage{pstricks}
\usepackage{amssymb}
\usepackage{amsmath}
\usepackage{indentfirst}
\usepackage{epsfig}
\usepackage{indentfirst}
\usepackage{amscd}
\usepackage{multirow}
\usepackage[english,spanish]{babel}

\title{ Representaciones Circulares de  Grafos  Simples  Conexos y  el  Rango M\1nimo Semidefinido  de   un  Delta  Grafo  }

\author{Pedro D\1az Navarro\thanks{Escuela de Matem\'atica, Universidad de Costa Rica}}
\usepackage{ graphics}

\newcommand{\mathsym}[1]{{}}
\newcommand{\unicode}[1]{{}}

\date{Junio  de  2018}
\newtheorem{defi}{Definici\'on}

\newtheorem{teor}{Teorema}

\newtheorem{result}{Resultado.}
\newtheorem{example}{Ejemplo.}

 \font\conj =
msbm10 at 12pt  

\def\text#1{\mathop{{\rm{#1}}}\nolimits}

\def\msr{\mathop{{\rm{msr}}}\nolimits}

\def\cc{\mathop{{\rm{cc}}}\nolimits}

\def\d{\mathop{{\rm{d}}}\nolimits}

\def\Re{{\mbox{\conj R}}}

\def\C{{\mbox{\conj C}}}

\def\dis{\displaystyle}
\def\1{\'{\i}}
\begin{document}
\maketitle

\begin{abstract}
En  este  art\1culo  se presenta  la  t\'ecnica para  representar  grafos  simples y conexos por medio de   encustamientos  circulares  y se  muestra  su aplicaci\'on en  el  c\'alculo de   representaciones ortogonales   de  grafos simples conexos    esenciales en el estudio  de   conjeturas  de la teor\1a del an\'alisis matricial  y  la  teor\1a  de   grafos.
\end{abstract}
\bigskip
\noindent{\bf Palabras  clave:} conjetura delta, grafo  simple conexo, rango m\1nimo semidefinido, $\delta$-grafo, C-$\delta$ grafo, representaci\'on  ortogonal.\\
\\

\noindent{\bf Clasificaci\'on  Matem\'atica:} 05C50,05C76 ,05C85 ,68R05 ,65F99,97K30.

\section{Introducci\'on}

En la teor\1a  de  grafos   se  diferentes  t\'ecnicas de tipo  combinatorio  para probar resultados.  Algunas de ellas se basan  en el  coloramiento   de  los v\' ertices y  el c\'alculo  de  par\'ametros  determinados por el  m\1nimo  obtenido en alg\'un proceso  combinatorio, de  cubrimientos  de un grafo por medio de  familias   de otros  grafos, por separaci\'on  en  subfamilias  de   grafos inducidos, etc.

Estas  t\'ecnicas  requieren  de una  gran cantidad de  c\'alculo  en  virtud de que se  debe analizar  exhaustivamente  el  grafo para lograr  obtener el   par\'ametro  deseado. Si el  grafo  tiene  un n\'umero de vertices peque\~no, el c\'alculo  se podr\1a realizar  con mayor o menor  dificultad  requiri\'endose en la mayor\1a de los  casos  una  gran cantidad de tiempo.  Sin  embargo,  cuando  se  trabaja  con  grafos con  un n\'umero de v\'ertices  muy  grande o  bien  con   familias  infinitas de  grafos,  estas  t\'ecnicas    representan  un  verdadero desaf\1o  para  el  investigador.  Un n\'umero muy  elevado de v\'ertices  implica  necesariamente   una  gran cantidad de tiempo de proceso y,  dependiendo  de lo complicado que sean los c\'alculos, obtener  resultados  ser\'a  una tarea  dificil  incluso usando  ordenadores de alto desempeño    y los   mejores programas  que se tengan a mano. M\'as a\'un, dado  que la  representaci\'on de un  grafo  no impone  ninguna restricci\'on  en especial  sobre   la  ubicaci\'on  de los  v\'ertices  ni  de la  longitud   o  forma  de sus aristas,   la forma  en que se represente el  grafo puede  complicar aun m\'as  su estudio.

Por lo tanto,  es necesario  representar  un  grafo de  forma  tal  que permita identificar     caracter\1sticas  como por ejemplo  conectividad,  simetr\1as, grafos inducidos,  caminos maximales etc.

En este  art\1culo se  expone  una  manera de representar un  grafo $G$   llamada  {\bf  encrustamientos  circulares en el  sentido  horario}. Esta  t\'ecnica    permite no  solo  visualizar    el  grafo en una  forma ordenada    sino  tambi\'en  permite obtener informaci\'on  de este, de sus  subgrafos inducidos  y de su  grafo complementario.

Adem\'as  se mostrar\'a   como esta  t\'ecnica permite  la demostraci\'on de la conjetura  delta    para algunas familias  infinitas de  grafos simples  y  conexos constituyendose as\1 en una herramienta  valiosa   en el  estudio de la  teor\1a  de  grafos.

\section{Preliminares  de la  teor\1a de  grafos}
\addtocontents{toc}{\vspace{-15pt}}
En esta  secci\'on   daremos  algunas definiciones  y  resultados    de la teor\1a  de  grafos las cuales ser\'an usadas  en las secciones siguientes. Mayores detalles  se pueden  encontrar en  \cite{BO,BM, CH}.

Un  {\bf grafo} {$G(V,E)$} es un par  ordenado  {$(V(G),E(G)),$} donde  {$V(G)$} es el conjunto de v\'ertces  y   {$E(G)$} es  el  conjunto de aristas  junto  con una  {\bf funci\'on de incidencia} $\psi(G)$ que asocia cada arista   de  $G$ un  par no ordenado de  v\'ertices (no necesariamente distintos) de  $G$. El {\bf orden} de {$G$}, denotado {$|G|$}, es el n\'umero de  v\'ertices en {$G$}.  Un  grafo se dice ser   {\bf simple} si no  tiene bucles o  m\'ultiples aristas  entre dos  v\'ertices  dados. El  {\bf complemento} de un grafo {$G(V,E)$} es el  grafo  {$\overline{G}=(V,\overline{E}),$} donde  {$\overline{E}$} consiste de todas las aristas  que no est\'an  en  {$E$}.
Un {\bf  subgrafo} {$H=(V(H),E(H))$} of {$G=(V,E)$} es un grafo con  {$V(H)\subseteq V(G)$} y {$E(H)\subseteq E(G)$}. Un {\bf subgrafo inducido} {$H$} de {$G$}, denotado G[V(H)], es un subgrafo  con  {$V(H)\subseteq V(G)$} and {$E(H)=\{\{i,j\} \in E(G):i,j\in V(H)\}$}. Algunas  veces denotamos la arista $\{i,j\}$ como  $ij$.
Decimos que  dos  v\'ertices   de un  grafo  $G$ son {\bf adyacentes}, denotado  $v_i\sim v_j$,   si  existe  una arista $\{v_i,v_j\}$ en $G$.  De otra  forma  decimos que dos  v\'ertices    $v_i$ and $v_j$ son  {\bf no adjacentes}  y  denotamos esto por   $v_i \not\sim v_j$.  Sea $N(v)$ el conjunto de  v\'ertices  que son adyacentes al  v\'ertice  $v$ y sea  $N[v]=\{v\}\cup N(v)$.

 El  {\bf grado } de un  v\'ertice  $v$ en $G$, denotado $\d_G(v)$, es la cardinalidad de  $N(v)$. Si  $\d_G(v)=1$  entonces  $v$ se llama  {\bf v\'ertice  colgante} de $G$. Usamos  $\delta(G)$ para denotar el  grado m\1nimo  de los  v\'ertices en   $G$,mientras que  $\Delta(G)$ denot\'a el m\'aximo grado de los v\'ertices en  $G$.

Dos  grafos $G(V,E)$ and $H(V',E')$  son  {\bf id\'enticos} , y  se denota ,  $G=H$, si  $V=V',  E=E'$, y $\psi_G=\psi_H$. Dos  grafos  $G(V,E)$ y $H(V',E')$ son {\bf isomorfos}, y  se denota por   $G\cong H$, si existen  biyecci\'ones  $\theta:V\to V'$  y  $\phi: E\to  E' $ tal  que   $\psi_G(e)=\{u,v\}$  si  y solo si   $\psi_H(\phi(e))= \{\theta(u), \theta(v)\}$.
 Un {\bf grafo  completo}  es un  grafo  simple  en el  cual  los  v\'ertices  son  dos a dos  adyacentes.
  Usaremos  {$nG$} para denotar  {$n$} copias de un  rafo {$G$}. Por ejemplo , $3K_1$  denota tres  v\'ertices aislados  $K_1$ mientras   que {$2K_2$} es el  grafo dado por  dos  copias  disconexas de  $K_2$.
 Un  {\bf camino } es  una lista   de  v\'ertices distintos  en la  cual  los  v\'ertices sucesivos estan  conectados por una arista. Un camino de   {$n$} v\'ertices se denota  {$P_n$}. Un  grafo {$G$} se dice ser  {\bf conexo} si hay un camino entre cualesquiera  dos  v\'ertices de  {$G$}. Un  {\bf ciclo} de {$n$} v\'ertices, denotado {$C_n,$} es un camino  donde el punto inicial  y  final son el mismo. Un   {\bf \'arbol} es un grafo  conectado sin ciclos. Un grafo  $G(V,E)$ se llama  {\bf cordal}  si no  tiene ciclos inducidos  $C_n$ con  $n\ge 4$.
  Una {\bf componente}  de un  grafo $G(V,E)$ es un  subgrafo maximal  conexo. Un {\bf v\'ertice de corte} es un v\'ertice cuya eliminaci\'on  incrementa el n\'umero de componentes.
la  {\bf uni\'on}  $G\cup G_2$ de dos  grafos  $G_1(V_1,E_1)$ y $G_2(V_2,G_2)$  es la uni\'onde sus  conjuntos  de v\'ertices  y aristas, esto es  $G\cup G_2(V_1\cup V_2,E_1\cup E_2$. Donde   $V_1$ y $V_2$ son  conjuntos  disjuntos  la uni\'on  se llama  {\bf  uni\'on disjunta } y se denota  $G_1\sqcup G_2$.


\section{El  Rango M\1nimo Semidefinido  de  un  Grafo}
En esta  secci\'on  establecemos   algunos  de los  resultados  conocidos para   e;  rango minimo semidefinido ($\msr$)de un  grafo  $G$  que ser\'a  usado en las secciones siguientes.

Una  {\bf matriz positiva}  $A$ es  una matriz  Hermitiana    $n\times n$ tal que  $x^\star A x>0$  para  todo  vector  $x\in \C^n$ no nulo. Equivalentemente,  $A$  es una matriz  Hermitiana  definida positiva   $n\times n$  si  y solo si  todos los   eigenvalores   de $A$ son  positivos (\cite{RC}, p.250).

Una matriz Hermitiana $n\times n$ tal  que   $x^\star A x\ge 0$ para  todo  $x\in \C^n$ se llama  {\bf  semidefinida positiva  (psd)}. Equi\-va\-lentemente,   $A$ es  una matriz  Hermitiana  Semidefinida Positiva  $n\times n$  si  y solo si   $A$  has  all eigenvalores  no negativos (\cite{RC}, p.182).

Si $\overrightarrow{V}=\{\overrightarrow{v_1},\overrightarrow{v_2},\dots, \overrightarrow{v_n}\}\subset \Re^m$  es  un  conjunto de vectores  columna  entonces la matriz
$ A^T A$, where $A= \left[\begin{array}{cccc}
  \overrightarrow{v_1} & \overrightarrow{v_2} &\dots& \overrightarrow{v_n}
\end{array}\right]$
y $A^T$  representa  la matriz  transpuesta  de   $A$, es una matriz  psd llamada la  {\bf matriz de Gram } of $\overrightarrow{V}$. Sea $G(V,E)$  un  grafo asociado  con una matriz  de Gram. Entonces  $V_G=\{v_1,\dots, v_n\}$ corresponde a un  conjunto de vectores  en  $\overrightarrow{V}$ and  E(G) corresponde a los productos  internos  no nulos  entre los  vectores  en    $\overrightarrow{V}$. En este caso $\overrightarrow{V}$  es llamada  una   {\bf  representaci\'on ortogonal } de  $G(V,E)$ in $\Re^m$. Si  tal  representaci\'on   ortogonal existe para     $G$ entonces  $\msr(G)\le m$.
Se  dice   que  dos  vectores $ \overrightarrow{v_i}, \overrightarrow{v_j}$   de  una  representaci\'on ortogonal    de un   grafo  simple  conexo  satisfacen una    {\bf relaci\'on  de ortogonalidad} si $ \overrightarrow{v_i}\cdot  \overrightarrow{v_j}=0$.  En  caso  contrario  se dice  que  satisfacen una {\bf relaci\'on  de adyacencia}.

Algunos de los  resultados  m\'as  comunes   acerca  del  rango m\1nimo  semidefinido   de  un  grafo   son los siguientes:

\begin{result}\cite{VH}\label{msrtree}
{Si  $T $ es un \'arbol entonces  $\msr(T)= |T|-1$.}
\end{result}
\begin{result}\cite{MP3}\label{msrcycle}
 {El  Ciclo $C_n$ tiene  rango m\1nimo semidefinido  $n-2$.}
\end{result}


\begin{result}\label{res2}
 \cite{MP3}{ \ Si un  grafo conexo  $G$  tiene  un  vertice  colgante  $v$, entonces  $\msr(G)=\msr(G-v)+1$ donde  $G-v$  es el    subgrafo  inducido  de  $G$ al eliminar  $v$.}
\end{result}

\begin{result} \cite{PB} \label{OS2}
 {Si {$G$} es un  grafo  cordal  y  conexo entonces $\msr(G)=\cc(G).$}
\end{result}

\begin{result}\label{res1}
\cite{MP2}\ {Si un  grafo  $G(V,E)$ tiene un  v\'ertice de corte tal  que , $G=G_1\cdot G_2$, entonces   $\msr(G)= \msr(G_1)+\msr(G_2)$.}
\end{result}

Las dos definiciones  siguientes nos dan  dos  familias  de  grafos   que son importantes  en el estudio  del  rango m\1nino semidefinido  de un  grafo  simple conexo.  La primera es  la definici\'on  de un   $\delta$-grafo.

\begin{defi}\label{ccpg}{Suponga  que  $G=(V,E)$  con  $|G|=n \ge 4$  es simple  y  conexo  tal que    $\overline{G}=(V,\overline{E})$  es tambi\'en simple  y  conexo. Decimos que  $G$ es un   {\bf $\mathbf{\delta}$-grafo} si  se  puede  etiquetar los  v\'ertices  de $G$ de forma tal  que:
\begin{enumerate}
  \item[(1)] el  grafo inducido de los  v\'ertices   $v_1,v_2,v_3$  en $G$ es  o bien  $3K_1$  o   $K_2 \sqcup K_1$, y
  \item[(2)] para  $m\ge 4$,  el v\'ertice  $v_m$  es adjacente a  todos los  v\'ertices  previos  $v_1,v_2,\dots,v_{m-1}$  excepto a  lo  sumo $\dis{\left\lfloor\frac{m}{2}-1\right\rfloor}$ v\'ertices.
 \end{enumerate}}
\end{defi}
 La segunda   se refiere  a una familia de  grafos  que contiene al  complemento de un  $\delta$-grafo.

\begin{defi} {Suponga  que    un grafo $G(V,E)$ con $|G|=n \ge 4$  es  simple  y  conexo   tal  que   $\overline{G}=(V,\overline{E})$ es tambi\'en  simple  y  conexo. Decimos que   $G(V,E)$  is a {\bf C-$\delta$-grafo} si    $\overline{G}$   es un  $\delta$-grafo.
en  otras palabras,  $G$   es un   {\bf C-$\mathbf{\delta}$ graph} si podemos etiquetar los  v\'ertices  de    $G$  de forma tal  que
 \begin{enumerate}
 \item[ (1)] el  grafo  inducido de los  v\'ertices   $v_1,v_2,v_3$ in  $G$ es o  bien   $K_3$ o $P_3$,  y
 \item[(2)] para   $m\ge 4$,  el v\'ertice  $v_m$  es  adyacente a  a lo sumo   $\dis{\left\lfloor\frac{m}{2}-1\right\rfloor}$ de los  v\'ertices previos.
\end{enumerate}}
\end{defi}

\begin{example}\label{examplecp}
El  ciclo  $C_n, n\ge 6$  es  un  C-$\delta$ grafo y  su  complemento  es un    $\delta$-grafo.
\begin{center}
\includegraphics[height=50mm]{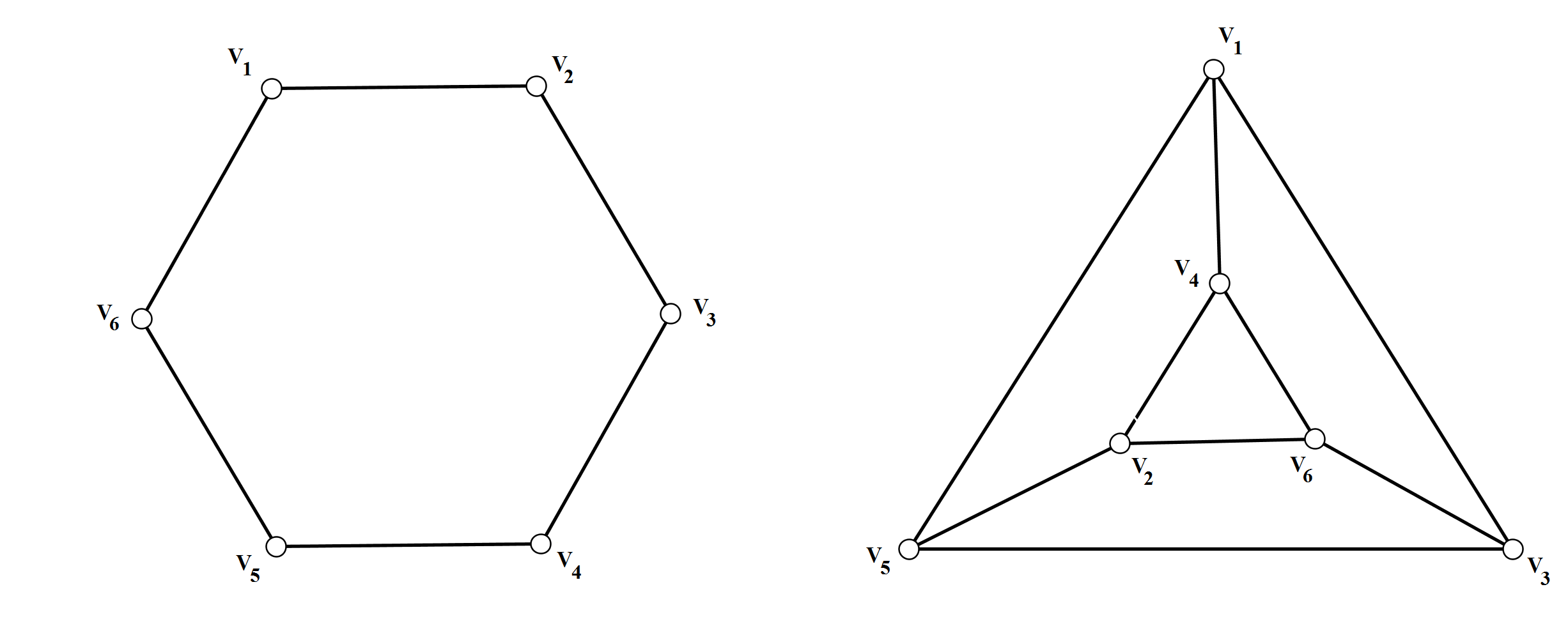}
 \end{center}
\vspace{-0.1in}\begin{figure}[h]
\centering
\caption{The  Graphs $C_6$  and  the  3-Prism }\label{fig3.1}
\end{figure}
Note que podemos etiquetar  los v\'ertices  de $C_6$ en  sentido  horario $v_1,v_2,v_3,v_4,v_5,v_6$. El  grafo inducido por  $v_1,v_2,v_3$  es $P_3$. El    v\'ertice $v_4$ es adyacente  al  v\'ertice previo  que es  $v_3$. Also,   v\'ertice $v_5$ es adyacente  a $v_4$ y  el v\'ertice  $v_6$ es adyacente  a dos de los  v\'ertices  previos  $v_1$ y $v_5$.  As\1, $C_6$  es un  C-$\delta$-graph. El $3$-prisma,  el  cual es isomorfo  a  el  complemento  de    $C_6$,  es un   $\delta$-grafo.
\end{example}

 En   \cite{PD} se prob\'o   el  siguiente  resultado
\begin{teor} \label{main} {Sea  $G(V,E)$ un  $\delta$-grafo entonces
$$
\msr(G)\le\Delta(\overline{G})+1=|G|-\delta(G)\label{mrsineq1}
$$}
\end{teor}

La   conjetura delta postulada   en  2006 en el     AIM workshop de  2006 in Palo Alto, California,  indica que  para  todo  grafo  $G$ simple  y  conexo  se tiene  que  $\msr(G)\le |G|-\delta(G)$. Detalles de esta  conjetura  se puede encontrar en  \cite{RB1}.
 Dado  que  el $\msr(G)$  es la m\1nima   dimensi\'on  en la que se puede encontrar una reprentaci\'on  ortogonal para el grafo  $G$,  el   teorema anterior  claramente establece la  validez   de esta conjetura para  todo  $\delta $-grafo.
 \section{Encrustamientos circulares  de  grafos }

 Como  se  indic\'o  anteriormente,  la  representaci\'on de  un  grafo simple conexo no  depende  de la posici\'on  en que se ubiquen los  v\'ertices  ni  de la  longitud  o  forma de  las aristas que unen   estos  v\'ertices   En esta secci\'on  mostramos   como representar   un   grafo por medio de  un encrustamiento de este  en  un c\1rculo dado.

 \begin{defi}\label{CCR}
Sea  $G(V_1,E_1)$ un  grafo    simple conexo. Sea  $V=\{v_1,v_2,\dots,v_n\}$  el  conjunto de  v\'ertices  de $G$. La  {\bf representaci\'on   circular  en sentido  horario} del  grafo  $G$  es la  representaci\'on  gr\'afica que se obtiene  cuando  sus   v\'ertices   $v_1,v_2,\dots,v_n$ , se colocan sobre un c\1rculo de forma tal  que dos  v\'ertices  consecutivos   subtiendan un  arco de $\frac{2\pi}{n}$ y  se colocan  en sentido  horario en  orden  creciente  de \1ndices.
\end{defi}

\begin{example}
{La  escalera de M\"obius  de orden  $2n$, $ML_{2n}$ tiene  tradicionalmente la  siguiente  representaci\'on la cual  es b\'asicamente  una escalera  incrustada en una cinta de m\"obius.\\
\begin{center}
\includegraphics[height=30mm]{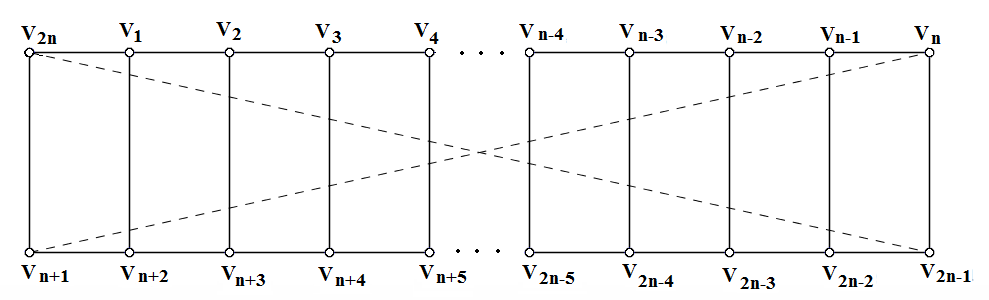}
\end{center}
\vspace{-0.1in}\begin{figure}[h]
\centering
\caption{Escalera de M\"obius de orden $2n$} \label{ML2nt}
\end{figure}
}
\end{example}
 Aplicando la  definici\'on  de la  representaci\'on  circular  en sentido  horario  obtenemos  la  siguiente  representaci\'on para $ML_{2n}$\\
 \begin{center}
\includegraphics[height=55mm]{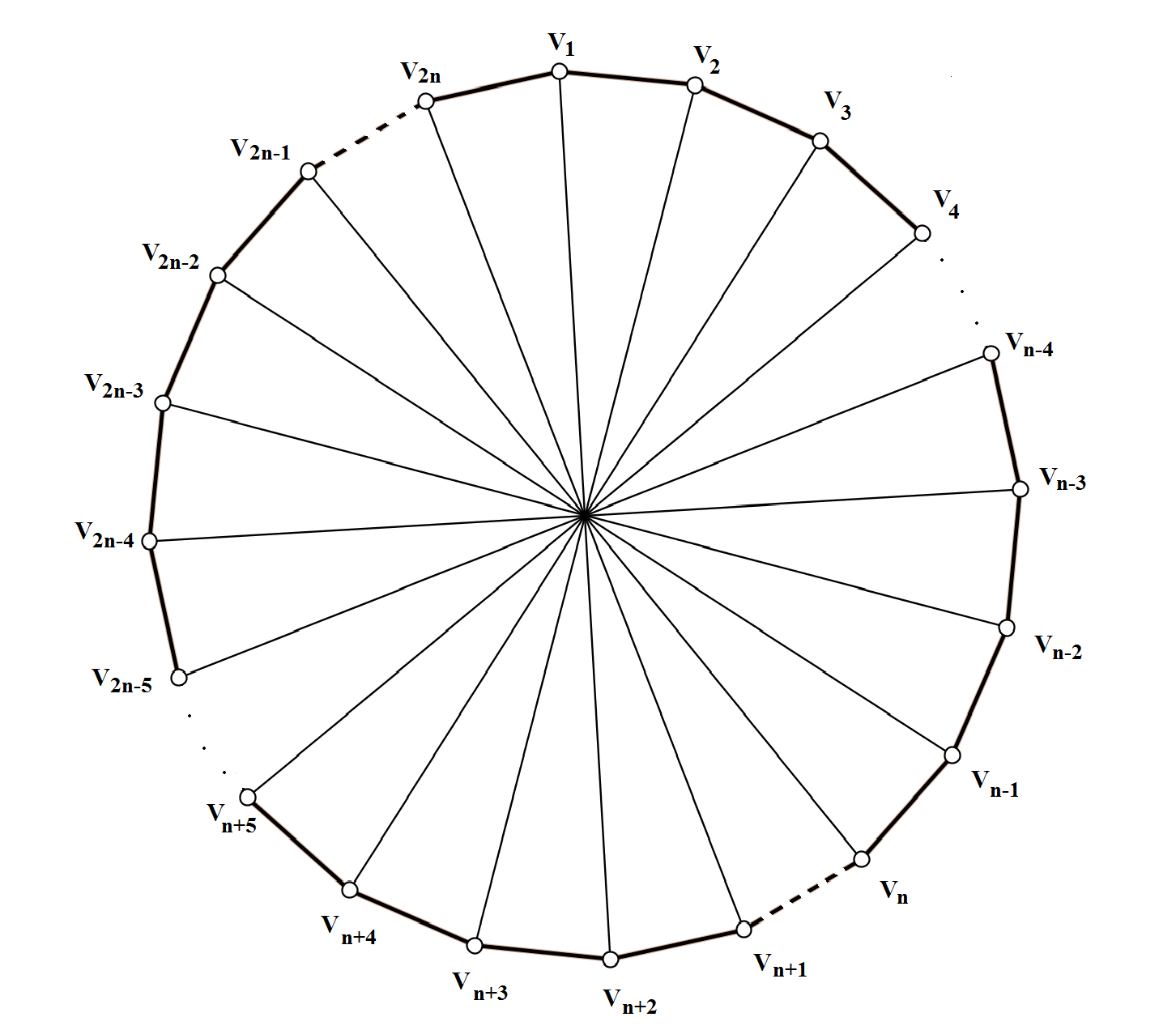}
\end{center}
\vspace{-0.1in}\begin{figure}[h]
\centering
\caption{Escalera de M\"obius de orden $2n$} \label{rchML2n}
\end{figure}

 Las  representaciones  circulares en sentido  horario  son  especialmente \'utiles  cuando  representamos operaciones  con  grafos.  Un ejemplo de esto es  el caso de  la operaci\'on  llamada {\bf producto  cartesiano}   que se  define    de la siguiente manera:
\begin{defi}\label{capo}
El   {\bf Producto Cartesiano} $G\Box H$  de dos  grafos simples  y  conexos  $G(V_1,E_1)$ y $H(V_2,E_2)$  es el  grafo  cuyo  conjunto de  v\'ertices es  $V_1\times V_2$ y  el  conjunto de aristas  se define  como sigue:
\begin{equation}
 (u,v)\sim(u',v')\Leftrightarrow
\left\{
  \begin{array}{ll}
   ( u=u') \  \wedge\ (vv'\in E(H)) &  \\
   ( v=v')\ \wedge\ (uu'\in E(G))&
  \end{array}
\right.
\end{equation}
\end{defi}
Como  caso particular  consideremos   el producto  cartesiano   de  un  grafo  completo  $K_3$  y  un  camino  $P_4$. Usando  una  forma  tradicional para  representar   este   grafo  y  con  un poco de  trabajo podr\1amos,  despu\'es de  cierto trabajo  obtener una  representaci\'on    como la siguiente  en  donde se aprecia  que  \'esta  operaci\'on  produce   cuatro copias   de  $K_3$ y tres copias  de $P_4$. Cada una de las cuatro  copias de  $K_3$    tiene un  v\'ertice   en  cada una  de las  copias  de $P_4$.
 \begin{center}
\includegraphics[height=35mm]{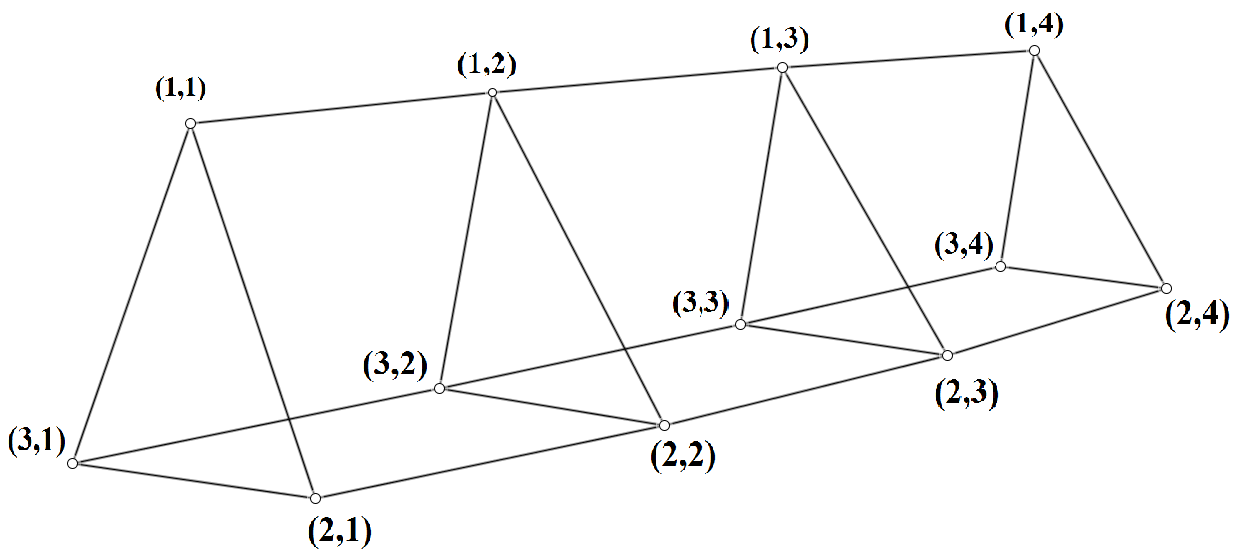}
\end{center}
\vspace{-0.1in}\begin{figure}[h]
\centering
\caption{Producto Cartesiano  $K_3\square P_4$} \label{k3cpp4t}
\end{figure}
Aplicando  la definici\'on de  la  representaci\'on   circular  horaria  a  $K_3\Box P_4$  se obtiene  el  siguiente  grafo equivalente   al  grafo de la  figura  \ref{k3cpp4t}.
\begin{center}
\includegraphics[height=45mm]{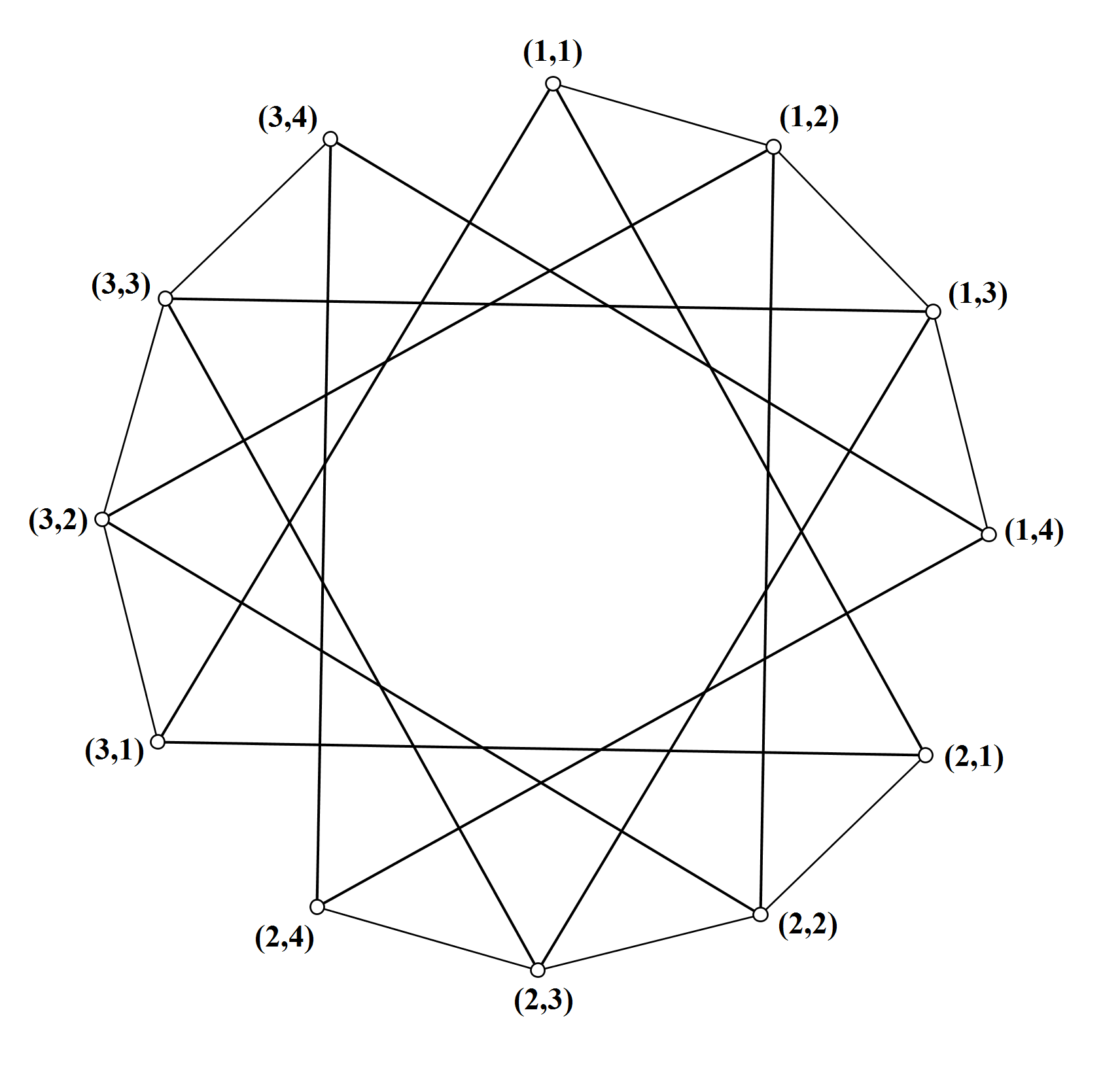}
\end{center}
\vspace{-0.1in}\begin{figure}[h]
\centering
\caption{Representaci\'on  circular  del Producto Cartesiano  $K_3\square P_4$ } \label{k3cpp4a}
\end{figure}
Mas  a\'un,   las   representaciones  circulares  horarias de  grafos  son    extremadamente  \'utiles   para  determinar  el  grafo  complementario $\overline{G}$     de  un  grafo  $G$  ya que  la posicion  de los  v\'ertices en el  c\1rculo es la misma para ambos  grafos  y  el  trazo  del  grafo complementario  se obtiene  f\'acilmente de observar las aristas  faltantes  del  grafo  completo  $K_n$  del  cual  $G$  es un  subgrafo  inducido. Por ejemplo,  tomando  de  nuevo   la  representaci\'on  del producto  cartesiano $K_3\Box P4$  de la  figura \ref{k3cpp4t}, si   colocamos los  v\'ertices de $K_3\Box P4$ en la misma posici\'on   y dibujamos   el  grafo  complementario  $\overline{K_3\Box P4}$   se obtiene  una representaci\'on  como la siguiente:
\begin{center}
\includegraphics[height=45mm]{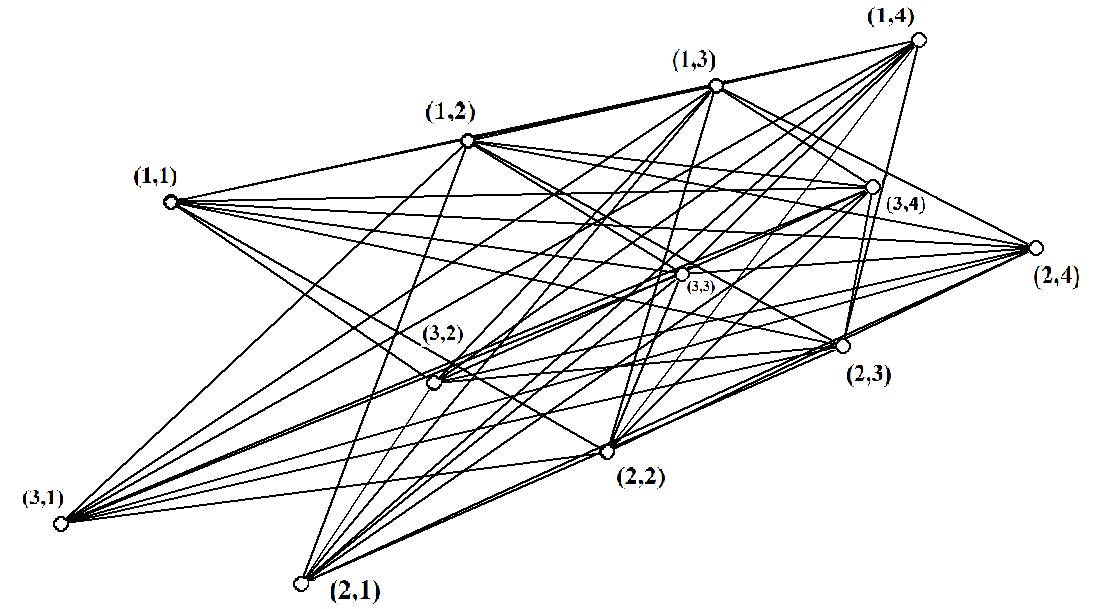}
\end{center}
\vspace{-0.3in}\begin{figure}[h]
\centering
\caption{Grafo Complementario  del Producto Cartesiano  $K_3\square P_4$} \label{cok3cpp4t}
\end{figure}
 Sin  embargo,  si  representamos  el  grafo $\overline{K_3\square P_4}$ usando la  representaci\'on  circular  en sentido  horario  se obtiene
\begin{center}
 \includegraphics[height=45mm]{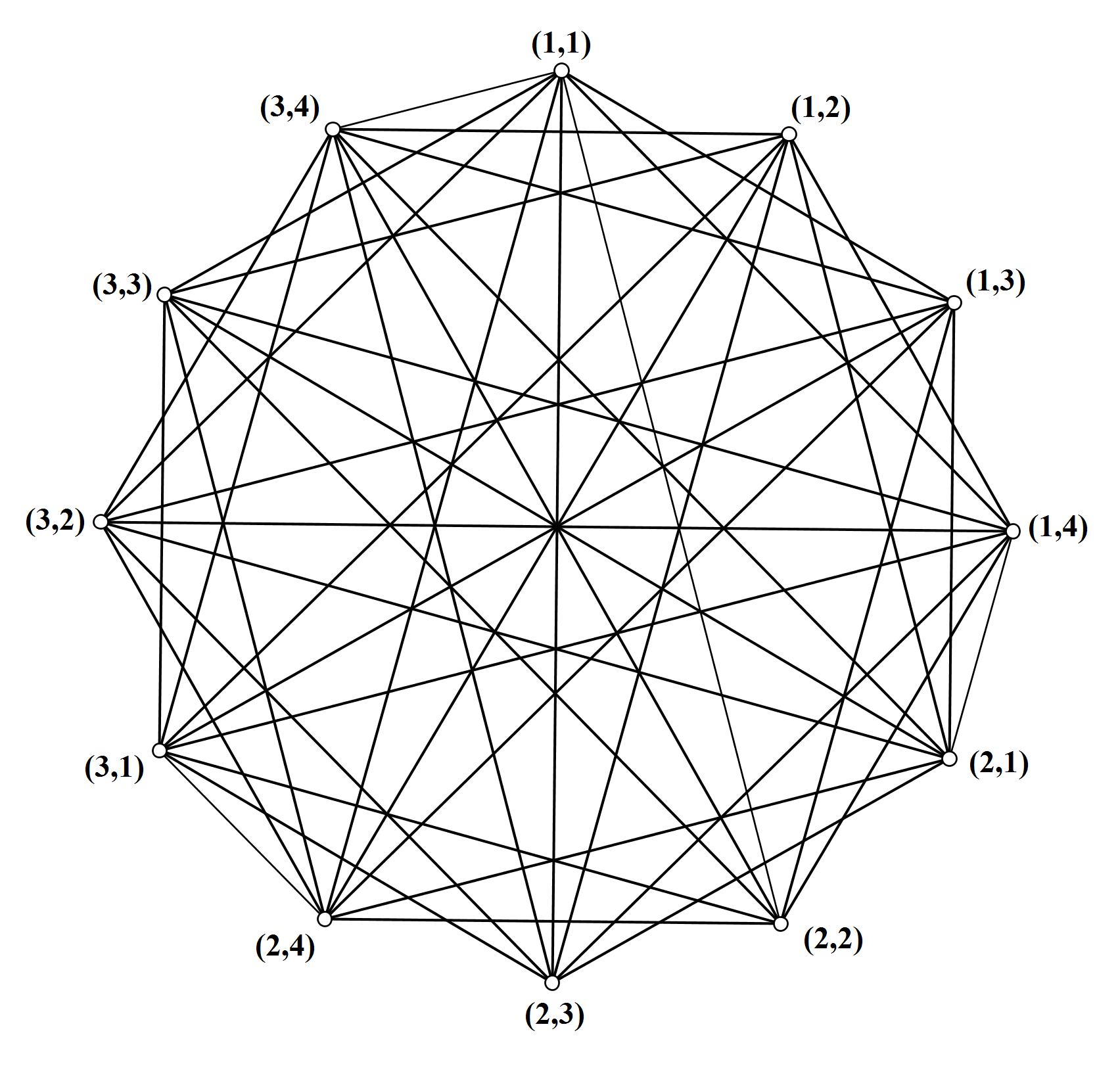}
\end{center}
\vspace{-0.3in}\begin{figure}[h]
\centering
\caption{Grafo Complementario  del Producto Cartesiano  $K_3\square P_4$} \label{cok3cpp4t}
\end{figure}
Claramente,  es  una  forma m\'as sencilla  de representar  $\overline{K_3\square P_4}$ desde la cual se pueden   observar  las  relaciones de  adyacencia    entre los  v\'ertices.

\section{Representaci\'ones circulares y  rango m\1nimo semidefinido }

Las  representaciones  circulares   de   grafos  simples conexos  son  utiles   en el   an\'alisis  matricial  y la  teor\1a de grafos pues nos pernite  calcular de manera  relativamente  sencilla una  representaci\'on  ortogonal   de cualquier  grafo simple conexo.  La  complejidad de los  c\'alculos  depende de lo  complejo  que sea  el  grafo  simple conexo.
Consideremos  por  ejemplo    el  grafo $K_3\square P_4$ y  su  grafo  complementario. Despu\'es de  realizar  el  encrustamiento  circular    de ambos  grafos  las  representaciones de estos   est\'a  dada en la  siguiente  figura.
\begin{center}
\end{center}
\begin{center}
\includegraphics[height=45mm]{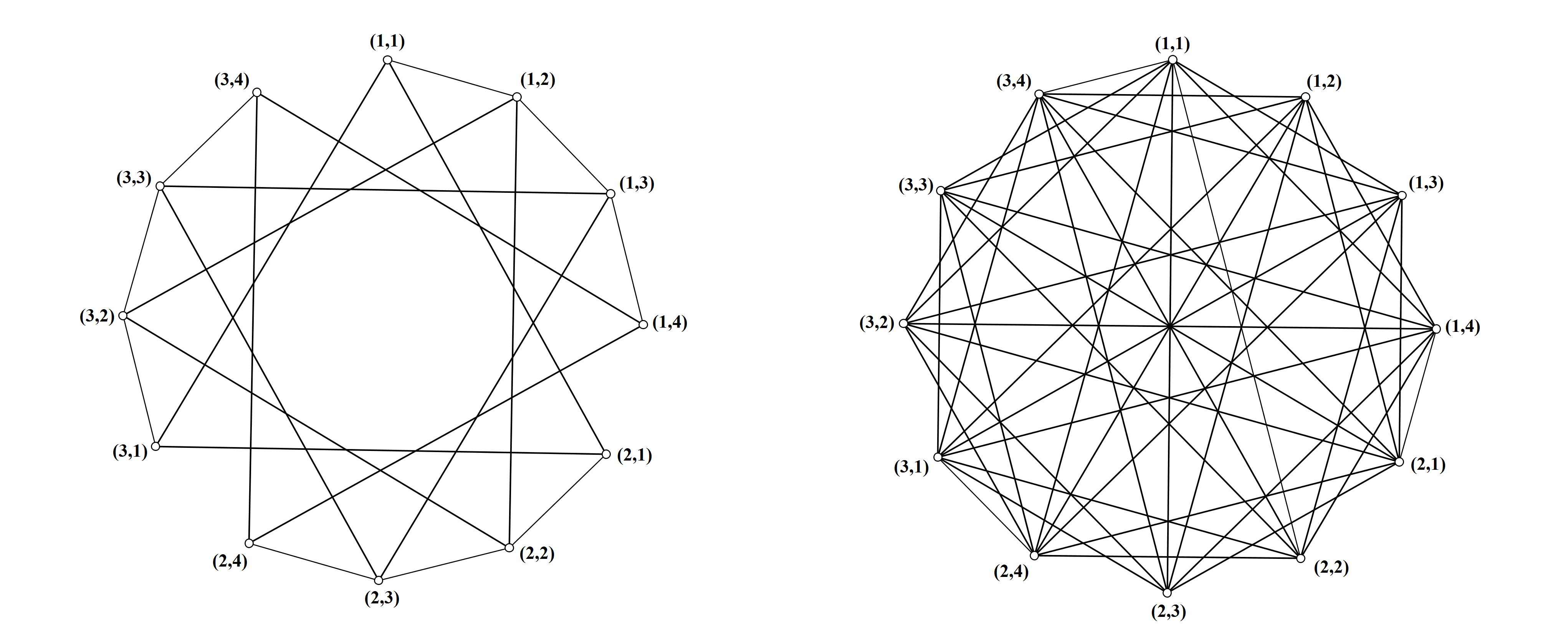}
 \end{center}
\vspace{-0.1in}\begin{figure}[h]
\centering
\caption{ Grafo  $K_3\square P_4$ su  complemento  $\overline{K_3\square P_4}$ }
\label{figB.1.1}
\end{figure}
 De \cite{PD} se  sabe  que   $\Delta(K_n\square P_4)= n+1$  y  adem\'as es   un $\delta$- grafo, por lo  tanto,  en  virtud  del  teorema  \ref{main} se  tiene   que    $\msr(K_3\square P_4)\le  5$.  En  consecuencia debe  ser posible encontrar   una  representaci\'on ortogonal  de  vectores  linealmente independientes  para el  grafo  $K_3\square P_4$.

 La  t\'ecnica descrita  en la  demostraci\'on  de   \ref{main}  en  \cite{PD}  para  calcular la  representaci\'on  ortogonal  de un  $\delta$- grafo consiste en  calcular las  representaciones  ortogonales de    cada uno de los   grafos  inducidos    de  una  secuencia  $Y_1\subseteq Y_2\subseteq\dots \subseteq Y_n= k_3\square P_4 $   de  $K_3\square P_4$ que se obtienen al  ir agregando  un  v\'ertice  a la  vez y  calcularle  su  respectiva  representaci\'on  vectorial de  forma  tal que se satisfagan las  condiciones de adyacencia  y ortogonalidad  del  subgrafo inducido. 
  
  Sean   $e_1, e_2,\dots,  e_5$  son  los vectores  can\'onicos  de  $\Re^5$. El primer paso  consiste en enumerar los  v\'ertices  de  de  $K_3\square P_4$  en la  forma  $v_1=(1,1),v_2=(1,2),\dots, v_12=(3,4)$ siguiendo el sentido  horario. Con  esto se  genera una  secuencia  de   subgrafos  inducidos  $Y_1(v_1), Y_2( v_1,v_2), \dots, Y_n(v_1,v_2,v_n)=K_3\square P_4$. Luego   aplicando  condiciones de  adyacencia  y ortogonalidad  definidas por el  grafo,  esto es  $v_i\cdot v_j=0$  si  $i\ne j$   y  $v_i\cdot v_j\ne 0$ se obtiene una secuencia de  sistemas  de ecuaciones  lineales  no homog\'eneos $A_i v_i=d_i,  i=1,2,\dots,n$  donde   algunas  de las  entradas del  vector $d_i,  i=1,\dots,n$   son no nulas.  Como no importa  cual es el  valor  no nulo  de las  entradas n o nulas de  $d_i$, estas se podr\1an  considerar  como  variables auxiliares  que permite  transformar el sistema no  homog\'eneo  en un  $A_i v_i=d_i,  i=1,2,\dots,n$  en un  sistema  homog\'eneo  $Aw=0,  w\in  \Re^{5+p}$  donde  $p$  es el n\'umero de    entradas no nulas  del  vector $d_i$  del   sistema  $A_i v_i=d_i,  i=1,2,\dots,n$. Esto  da  un margen de libertad  mayor para escoger las  soluciones   de las  entradas  del vector  $v_i$ y permite  encontrar  una  representaci\'on  ortogonal  de  $Y_i,  i = 1,2,\dots,n$  de  vectores  linealmente independientes.
  
  Si  suponemos  que  $\overrightarrow{v}_i= k_{i,1}\overrightarrow{e}_1 + k_{i,2}\overrightarrow{e}_2 + k_{i,3}\overrightarrow{e}_3 + k_{i,4}\overrightarrow{e}_4 + k_{i,5}\overrightarrow{e}_5$  es  la  representaci\'on  vectorial   del  vector  $\overrightarrow{v}_i,  i = 1,2,\dots,12$  de la  representac\'on  ortogonal  de  $K_3\square P_4$  que  buscamos    entonces el c\'alculo de la  representaci\'on ortogonal   de  $k_3\square P_4$  usando la  t\'ecnica  descrita  en  \cite{PD}    es el  siguiente:
  
  \begin{itemize}
\item {\bf V\'ertice $\mathbf v_1= (1,1)$}.\\
 Sea $\overrightarrow{v}_1 = k_{1,1}\overrightarrow{e}_1 + k_{1,2}\overrightarrow{e}_2 + k_{1,3}\overrightarrow{e}_3 + k_{1,4}\overrightarrow{e}_4 + k_{1,5}\overrightarrow{e}_5$. Podemos  comenzar  con  cualquier   vector no nulo.  Tome 
$$
\overrightarrow{v}_1=\overrightarrow{e}_1.
$$
\item {\bf V\'ertice $\mathbf v_2= (1,2)$}.\\
sea  $\overrightarrow{v}_2 = k_{2,1}\overrightarrow{e}_1 + k_{2,2}\overrightarrow{e}_2 + k_{2,3}\overrightarrow{e}_3 + k_{2,4}\overrightarrow{e}_4 + k_{2,5}\overrightarrow{e}_5$.
Como  $v_2$ no es adyacente  a  $v_1$ in $\overline{K_3\square P_4}$, se obtiene  una  sola ecuaci\'on  $k_{2,1} = 0$. todas las otras  entradas  son   variables libres. Entonces 
we choose
$$
\overrightarrow{v}_2 = \overrightarrow{e}_2.
$$
Claramente  $\langle v_1,v_2\rangle=0$ y se satisfacen las  condiciones de  adyacencia  y  ortogonalidad  en el  subgrafo inducido  $Y_2( V_1, v_2)$..
\item {\bf V\'ertice $\mathbf v_3= (1,3)$}.\\
  Sea $\overrightarrow{v}_3 = k_{3,1}\overrightarrow{e}_1 + k_{3,2}\overrightarrow{e}_2 + k_{3,3}\overrightarrow{e}_3 + k_{3,4}\overrightarrow{e}_4 + k_{3,5}\overrightarrow{e}_5$.

Como $v_3$ no es adyacente   a  $v_2$ en  $\overline{K_3\square P_4}$ pero es  adyacente  $v_1$ obtenemos las ecuaciones  $k_{3,1} = g_{3,1}, g_{3,1}\ne 0$
y   $k_{3,2} = 0$. Las dem\'as entradas  son  variables  libres . as\1  se obtiene la  siguiente matriz:
$$
A=\left(
\begin{array}{cccccc}
 1 & 0 & 0 & 0 & 0 & -1 \\
 0 & 1 & 0 & 0 & 0 & 0
\end{array}
\right)
$$
Entonces,  podemos escoger  $g_{3,1} = 1$. Luego $k_{3,1} = 1$ y $k_{3,2} = 0$. las otras dos entradas s on variables libre. Como el vector $v_3$ debe ser  dos a  dos  linealmente independiente   con  $v_1$ yand $v_2$ necesitamos escoger  otra nueva  variable libre  no nula. Sea $k_{3,3}=1$ entonces
$$
\overrightarrow{v}_3=\overrightarrow{e}_1+\overrightarrow{e}_3.
$$
As\1, $\langle \overrightarrow{v}_1,\overrightarrow{v}_3\rangle=1$ and $\langle \overrightarrow{v}_2,\overrightarrow{v}_3\rangle=0$  y se satisfacen las  condiciones de  adyacencia  y  ortogonalidad  en el  subgrafo inducido  $Y_3( V_1, v_2,v_3)$.
\item {\bf V\'ertice $\mathbf v_4= (1,4)$}.\\
Sea $\overrightarrow{v}_4 = k_{4,1}\overrightarrow{e}_1 + k_{4,2}\overrightarrow{e}_2 + k_{4,3}\overrightarrow{e}_3 + k_{4,4}\overrightarrow{e}_4 + k_{4,5}\overrightarrow{e}_5$.
Como $v_4$ no es adyacente a   $v_3$ en  $\overline{K_3\square P_4}$ pero es adyacente  a  $v_1$ and $v_2$ obtenemos las  ecuaciones  $k_{4,1} = g_{4,1}, k_{4,2} = g_{4,2}, g_{4,1}\ne 0, g_{4,2}\ne 0$   y  $k_{4,3} = 0$. De esta  forma obtenemos la matriz:
$$
\left(
\begin{array}{ccccccc}
 1 & 0 & 0 & 0 & 0 & -1 & 0 \\
 0 & 1 & 0 & 0 & 0 & 0 & -1 \\
 1 & 0 & 1 & 0 & 0 & 0 & 0
\end{array}
\right)\sim
\left(
\begin{array}{ccccccc}
 1 & 0 & 0 & 0 & 0 & -1 & 0 \\
 0 & 1 & 0 & 0 & 0 & 0 & -1 \\
 0 & 0 & 1 & 0 & 0 & 1 & 0
\end{array}
\right)
$$
De este  c\'alculo se puede escoger  $g_{4,1} = g_{4,2} = 1$ y$k_{4,3} =-1$. entonces $k_{4,1}=1, k_{4,2}=1$, y $k_{4,3}=-1$. Pero  en  tal caso   se hace imposible  satisfacer  la  condici\'on de ortogonalidad $\langle v_3, v_4 \rangle= 0$.  Como  los  sistemas para    $v_1$ and $v_3$ tienen  variables libres  podemos   usar una t\'ecnica  de rastreo,  redefinir   el  vector $\overrightarrow{v}_1$ as $\overrightarrow{v}_1= \overrightarrow{e}_1 +\overrightarrow{e}_4$   y  obtener la matriz :
$$
\left(
\begin{array}{ccccccc}
 1 & 0 & 0 & 1 & 0 & -1 & 0 \\
 0 & 1 & 0 & 0 & 0 & 0 & -1 \\
 1 & 0 & 1 & 0 & 0 & 0 & 0
\end{array}
\right)\sim
\left(
\begin{array}{ccccccc}
 1 & 0 & 0 & 1 & 0 & -1 & 0 \\
 0 & 1 & 0 & 0 & 0 & 0 & -1 \\
 0 & 0 & 1 & -1 & 0 & 1 & 0
\end{array}
\right)
$$
As\1, podemos  resolver  todas  las  condiciones de  ortogonalidad   y  adjacencia  para   obtener, $ k_{4,1} = -k_{4,4} + g_{4,1}, k_{4,2} = g_{4,2}, k_{4,3} = k_{4,4} - g_{4,1}$. Escogiendo 
$ g_{4,1} = 2, g_{4,2} = 1$, y $ k_{4,4}=1$ obtenemos  $k_{4,1}=1, k_{4,2}=1$, y $k_{4,3}=-1$.As\1
$$
\overrightarrow{v}_4= \overrightarrow{e}_1+\overrightarrow{e}_2-\overrightarrow{e}_3+\overrightarrow{e}_4.
$$
Entonces $ \langle v_1,v_4\rangle=2, \langle v_2,v_4\rangle=1, \langle v_3,v_4\rangle=0$.
En consecuencia, todas las condiciones  de  adyacencia  y  ortogonalidad   se  satisfacen.
\item {\bf V\'ertice $\mathbf v_5= (2,1)$}.\\
 Let $\overrightarrow{v}_5 = k_{5,1}\overrightarrow{e}_1 + k_{5,2}\overrightarrow{e}_2 + k_{5,3}\overrightarrow{e}_3 + k_{5,4}\overrightarrow{e}_4 + k_{5,5}\overrightarrow{e}_5$.

 El v\'ertice $v_5$ es adyacente a $v_2, v_3$, y a  $v_4$ en  $\overline{K_3\square P_4}$ pero no es  adyacente  a  $v_1$.  De estas  condiciones   se obtienen las  ecuaciones  $k_{5,1}+k_{5,4}=0, k_{5,2}-g_{5,2}=0
, k_{5,1}+k_{5,3}-g_{5,3}=0$, and $k_{5,1}+k_{5,2}-k_{5,3}+k_{5,4}-g_{5,4}=0$. Asi  obtenemos la  matriz:

$$
\left(
\begin{array}{cccccccc}
 1 & 0 & 0 & 1 & 0 & 0 & 0 & 0 \\
 0 & 1 & 0 & 0 & 0 & -1 & 0 & 0 \\
 1 & 0 & 1 & 0 & 0 & 0 & -1 & 0 \\
 1 & 1 & -1 & 1 & 0 & 0 & 0 & -1
\end{array}
\right)\sim
\left(
\begin{array}{cccccccc}
 1 & 0 & 0 & 0 & 0 & 1 & -1 & -1 \\
 0 & 1 & 0 & 0 & 0 & -1 & 0 & 0 \\
 0 & 0 & 1 & 0 & 0 & -1 & 0 & 1 \\
 0 & 0 & 0 & 1 & 0 & -1 & 1 & 1
\end{array}
\right).
$$
Entonces,  escogiendo  $g_{5,4} = 1, g_{5,3} = -1$ and $g_{5,2} = 2$ obtenemos $k_{5,4} = 2, k_{5,3} = 1, k_{5,2} = 2$, and $k_{5,1} = -2$. En  consecuencia, 
$$
\overrightarrow{v}_5=-2\overrightarrow{e}_1+2\overrightarrow{e}_2+\overrightarrow{e}_3+2\overrightarrow{e}_4.
$$
Luego $ \langle \overrightarrow{v}_5,\overrightarrow{v}_1\rangle=0,\langle \overrightarrow{v}_5,\overrightarrow{v}_2\rangle=2,\langle\overrightarrow{v}_5,\overrightarrow{v}_3\rangle=-1,\langle\overrightarrow{v}_5,\overrightarrow{v}_4\rangle=1
$  y las  condicines de adyacencia  y ortogonalidad se  satisfacen.
\item {\bf V\'ertice $\mathbf v_6= (2,2)$}.\\
Sea $\overrightarrow{v}_6 = k_{6,1}\overrightarrow{e}_1 + k_{6,2}\overrightarrow{e}_2 + k_{6,3}\overrightarrow{e}_3 + k_{6,4}\overrightarrow{e}_4 + k_{6,5}\overrightarrow{e}_5$.
El  v\'ertice $v_6$ es  adyacente  con $v_1, v_3$, y  $v_4$ in $\overline{K_3\square P_4}$ pero no es adyacente   con   $v_2$ and $v_5$. Deestas  condiciones  obten3emos las  ecuaciones $k_{6,1} + k_{6,4} -g_{6,1}=0, k_{6,2} = 0 , k_{6,1} + k_{6,3} - g_{6,3} = 0 ,k_{6,1} + k_{6,2} - k_{6,3} + k_{6,4} - g_{6,4} = 0$, and $-2k_{6,1}+2 k_{6,2}+k_{6,3}+2k_{6,4}=0$. La  matriz  de este  sistema es:
$$\left(
\begin{array}{cccccccc}
 1 & 0 & 0 & 1 & 0 & -1 & 0 & 0 \\
 0 & 1 & 0 & 0 & 0 & 0 & 0 & 0 \\
 1 & 0 & 1 & 0 & 0 & 0 & -1 & 0 \\
 1 & 1 & -1 & 1 & 0 & 0 & 0 & -1 \\
 -2 & 2 & 1 & 2 & 0 & 0 & 0 & 0
\end{array}
\right)\sim
\left(
\begin{array}{cccccccc}
 1 & 0 & 0 & 0 & 0 & 0 & -\frac{3}{7} & -\frac{2}{7} \\
 0 & 1 & 0 & 0 & 0 & 0 & 0 & 0 \\
 0 & 0 & 1 & 0 & 0 & 0 & -\frac{4}{7} & \frac{2}{7} \\
 0 & 0 & 0 & 1 & 0 & 0 & -\frac{1}{7} & -\frac{3}{7} \\
 0 & 0 & 0 & 0 & 0 & 1 & -\frac{4}{7} & -\frac{5}{7}
\end{array}
\right)
$$
As\1,   si   $g_{6,3} = g_{6,4} = 7$ we get $g_{6,1} = 9$ todas las condicioens de adyacencia  y ortogonalidad   se satisfacen. Tambi\'en , $k_{6,4} = 4, k_{6,3} = 2, k_{6,2} = 0$, and $k_{6,1}
= 5$. De esta  forma ,
$$
\overrightarrow{v}_6= 5\overrightarrow{e}_1+2\overrightarrow{e}_3+4\overrightarrow{e}_4.
$$
Luego,
$\langle v_6,v1\rangle=9, \langle v_6,v_2\rangle=0, \langle v_6,v_3\rangle=7, \langle v_6,v_4\rangle=7, \langle v_6,v_5\rangle=0$.
\item {\bf V\'ertice $\mathbf v_7= (2,3)$}.\\
Sea $\overrightarrow{v}_7 = k_{7,1}\overrightarrow{e}_1 + k_{7,2}\overrightarrow{e}_2 + k_{7,3}\overrightarrow{e}_3 + k_{7,4}\overrightarrow{e}_4 + k_{7,5}\overrightarrow{e}_5$.

El  v\'ertice $v_7$ es adyacente  con  $v_1, v_2, v_4$, y $v_5$ en  $\overline{K_3\square P_4}$ pero no es adyacente  a  $v_3$ y $v_6$. De estas  condiciones   obtenemos las  ecuaciones 
the equations: $k_{7,1}+k_{7,4}-g_{7,1}=0, k_{7,2}-g_{7,2}=0, k_{7_1}+k_{7,3}=0, k_{7,1}+k_{7,2}-k_{7,3}+k_{7,4}-g_{7,4}=0,-2k_{7,1}+2k_{7,2}+k_{7,3}+2k_{7,4}-g_{7,5}=0$, y $5k_{7,1}+2k_{7,3}+4k_{7,4}=0$. As\1 obtenemos la matriz: \\
$$
\left(
\begin{array}{ccccccccc}
 1 & 0 & 0 & 1 & 0 & -1 & 0 & 0 & 0 \\
 0 & 1 & 0 & 0 & 0 & 0 & -1 & 0 & 0 \\
 1 & 0 & 1 & 0 & 0 & 0 & 0 & 0 & 0 \\
 1 & 1 & -1 & 1 & 0 & 0 & 0 & -1 & 0 \\
 -2 & 2 & 1 & 2 & 0 & 0 & 0 & 0 & -1 \\
 5 & 0 & 2 & 4 & 0 & 0 & 0 & 0 & 0
\end{array}
\right)\sim
\left(
\begin{array}{ccccccccc}
 1 & 0 & 0 & 0 & 0 & 0 & 0 & -\frac{2}{7} & \frac{1}{7} \\
 0 & 1 & 0 & 0 & 0 & 0 & 0 & -\frac{9}{14} & -\frac{5}{28} \\
 0 & 0 & 1 & 0 & 0 & 0 & 0 & \frac{2}{7} & -\frac{1}{7} \\
 0 & 0 & 0 & 1 & 0 & 0 & 0 & \frac{3}{14} & -\frac{3}{28} \\
 0 & 0 & 0 & 0 & 0 & 1 & 0 & -\frac{1}{14} & \frac{1}{28} \\
 0 & 0 & 0 & 0 & 0 & 0 & 1 & -\frac{9}{14} & -\frac{5}{28}
\end{array}
\right)
$$
Entonces, si $g_{7,4} = 14, g_{7,5} =- 28$ luego $g_{7,2} = 4 , g_{7,1} = 2$   y  todas las  condiciones    de  adyacencia  y   ortogonalidad   se  satisfacen. Tambi\'en  $k_{7,4} = -6, k_{7,3} = -8, k_{7,2} = 4$, y $k_{7,1}=8$. As\1,
$$
\overrightarrow{v}_7=8\overrightarrow{e}_1+4\overrightarrow{e}_2-8\overrightarrow{e}_3-6\overrightarrow{e}_4.
$$
Entonces, $\langle \overrightarrow{v}_7,\overrightarrow{v}_1\rangle=2, \langle \overrightarrow{v}_7,\overrightarrow{v}_2\rangle=4,\langle \overrightarrow{v}_7,\overrightarrow{v}_3\rangle=0,\langle\overrightarrow{v}_7,\overrightarrow{v}_4\rangle=14,\langle \overrightarrow{v}_7,\overrightarrow{v}_5\rangle=-28,\langle \overrightarrow{v}_7,\overrightarrow{v}_6\rangle=0$.
\item {\bf V\'ertice $\mathbf v_8= (2,4)$}.\\
Sea $\overrightarrow{v}_8 = k_{8,1}\overrightarrow{e}_1 + k_{8,2}\overrightarrow{e}_2 + k_{8,3}\overrightarrow{e}_3 + k_{8,4}\overrightarrow{e}_4 + k_{8,5}\overrightarrow{e}_5$.

El  v\'ertice $v_8$ es adyacente  con   $v_1, v_2, v_3, v_5$, y  $v_6$ en  $\overline{K_3\square P_4}$ pero no es adyacente  con   $v_4$ y   $v_7$.  Aplicando estas  condiciones se obtienen las ecuaciones : $k_{8,1} + k_{8,4} - g_1 = 0, k_{8,2} - g_{8,2} = 0, k_{8,1} + k_{8,3}-g_{8,3} = 0, k_{8,1} + k_{8,2} - k_{8,3} + k_{8,4} = 0, -2 k_{8,1} + 2 k_{8,2} + k_{8,3} + 2 k_{8,4} - g_{8,5} = 0, 5 k_{8,1}+2 k_{8,3} + 4 k_{2,4} -g_{8,6}=0$, y $8k_{8,1}+4k_{8,2}-8k_{8,3}-6k_{8,4}=0$. La matriz  de este sistema es:
{\small $$
\left(
\begin{array}{cccccccccc}
 1 & 0 & 0 & 1 & 0 & -1 & 0 & 0 & 0 & 0 \\
 0 & 1 & 0 & 0 & 0 & 0 & -1 & 0 & 0 & 0 \\
 1 & 0 & 1 & 0 & 0 & 0 & 0 & -1 & 0 & 0 \\
 1 & 1 & -1 & 1 & 0 & 0 & 0 & 0 & 0 & 0 \\
 -2 & 2 & 1 & 2 & 0 & 0 & 0 & 0 & -1 & 0 \\
 5 & 0 & 2 & 4 & 0 & 0 & 0 & 0 & 0 & -1 \\
 8 & 4 & -8 & -6 & 0 & 0 & 0 & 0 & 0 & 0
\end{array}
\right)\sim
\left(
\begin{array}{cccccccccc}
 1 & 0 & 0 & 0 & 0 & 0 & 0 & 0 & \frac{2}{107} & -\frac{15}{107} \\
 0 & 1 & 0 & 0 & 0 & 0 & 0 & 0 & -\frac{49}{107} & -\frac{7}{107} \\
 0 & 0 & 1 & 0 & 0 & 0 & 0 & 0 & -\frac{33}{107} & -\frac{20}{107} \\
 0 & 0 & 0 & 1 & 0 & 0 & 0 & 0 & \frac{14}{107} & \frac{2}{107} \\
 0 & 0 & 0 & 0 & 0 & 1 & 0 & 0 & \frac{16}{107} & -\frac{13}{107} \\
 0 & 0 & 0 & 0 & 0 & 0 & 1 & 0 & -\frac{49}{107} & -\frac{7}{107} \\
 0 & 0 & 0 & 0 & 0 & 0 & 0 & 1 & -\frac{31}{107} & -\frac{35}{107}
\end{array}
\right)
$$}
As\1, si  $g_{8,5} = g_{8,6} = 107$ entonces $g_{8,3}=70, g_{8,2} = 56 , g_{8,1} = 3$  y   se satisfacen las condiciones de adyacencia  y ortogonalidad. Tambi\'en  $k_{8,4} =-16, k_{8,3} = 53, k_{8,2} =56$, and $k_{8,1} = 13$. Luego,
$$
\overrightarrow{v}_8=13 \overrightarrow{e}_1+56 \overrightarrow{e}2+53\overrightarrow{e}_3-16\overrightarrow{e}_4.
$$
En consecuencia, $\langle v_8, v_1\rangle=-3,\langle v_8,v_2\rangle=56,\langle v_8,v_3\rangle=66,\langle v_8,v_4\rangle=0, \langle v_8,v_5\rangle=107,\langle v_8,v_6\rangle=107,\langle v_8, v_7\rangle=0$.
\item {\bf V\'ertice $\mathbf v_9= (3,1)$}.\\
 Let $\overrightarrow{v}_9 = k_{9,1}\overrightarrow{e}_1 + k_{9,2}\overrightarrow{e}_2 + k_{9,3}\overrightarrow{e}_3 + k_{9,4}\overrightarrow{e}_4 + k_{9,5}\overrightarrow{e}_5$.

El v\'ertice $v_9$ es adyacente  con  $v_2, v_3, v_4, v_6,v_7$ y $v_8$ en $\overline{K_3\square P_4}$ pero no es adyacente  con  $v_1$ y $v_5$.  aplicando estas  condiciones  obtenemos las ecuaciones: $k_{9,1} + k_{9,4} = 0, k_{9,2} - g_2 = 0, k_{9,1} + k_{9,3} - g_{9,3} = 0, k_{9,1} + k_{9,2} - k_{9,3} + k_{9,4}-g_{9,4} = 0, -2 k_{9,1} + 2 k_{9,2} + k_{9,3} + 2 k_{9,4} = 0,  5 k_{9,1}+ 2 k_{9,3} + 4 k_{9,4}- g_{9,6}=0,  8 k_{9,1} + 4 k_{9,2} - 8 k_{9,3} - 6 k_{9,4}-g_{9,7} = 0$ y $13k_{9,1}+56k_{9,2}+53k_{9,3}-16k_{9,4}-g_{9,8}=0$. La matriz de este sistema  es:
{\tiny
$$
\left(
\begin{array}{ccccccccccc}
 1 & 0 & 0 & 1 & 0 & 0 & 0 & 0 & 0 & 0 & 0 \\
 0 & 1 & 0 & 0 & 0 & -1 & 0 & 0 & 0 & 0 & 0 \\
 1 & 0 & 1 & 0 & 0 & 0 & -1 & 0 & 0 & 0 & 0 \\
 1 & 1 & -1 & 1 & 0 & 0 & 0 & -1 & 0 & 0 & 0 \\
 -2 & 2 & 1 & 2 & 0 & 0 & 0 & 0 & 0 & 0 & 0 \\
 5 & 0 & 2 & 4 & 0 & 0 & 0 & 0 & -1 & 0 & 0 \\
 8 & 4 & -8 & -6 & 0 & 0 & 0 & 0 & 0 & -1 & 0 \\
 13 & 56 & 53 & -16 & 0 & 0 & 0 & 0 & 0 & 0 & -1
\end{array}
\right)\sim
\left(
\begin{array}{ccccccccccc}
 1 & 0 & 0 & 0 & 0 & 0 & 0 & 0 & 0 & -\frac{5}{392} & -\frac{1}{196} \\
 0 & 1 & 0 & 0 & 0 & 0 & 0 & 0 & 0 & -\frac{241}{3920} & -\frac{9}{1960} \\
 0 & 0 & 1 & 0 & 0 & 0 & 0 & 0 & 0 & \frac{141}{1960} & -\frac{11}{980} \\
 0 & 0 & 0 & 1 & 0 & 0 & 0 & 0 & 0 & \frac{5}{392} & \frac{1}{196} \\
 0 & 0 & 0 & 0 & 0 & 1 & 0 & 0 & 0 & -\frac{241}{3920} & -\frac{9}{1960} \\
 0 & 0 & 0 & 0 & 0 & 0 & 1 & 0 & 0 & \frac{29}{490} & -\frac{4}{245} \\
 0 & 0 & 0 & 0 & 0 & 0 & 0 & 1 & 0 & -\frac{523}{3920} & \frac{13}{1960} \\
 0 & 0 & 0 & 0 & 0 & 0 & 0 & 0 & 1 & \frac{257}{1960} & -\frac{27}{980}
\end{array}
\right)
$$
}
Entonces, si  $g_{9,8} = 1960, g_{9,7} = 3920$ then $g_{9,6} =-460, g_{9,4}=510, g_{9,3}=-200,g_{9,2} =250$ y  todas las condicines de adyacencia  y ortogonalidad se satisfacen. Tambi\'en $k_{9,4} = -60, k_{9,3}= -260, k_{9,2} =250$, y $k_{9,1} = 60$. en consecuencia
$$
\overrightarrow{v}_9=60\overrightarrow{e}_1+250\overrightarrow{e}_2-260\overrightarrow{e}_3-60\overrightarrow{e}_4.
$$
Entonces,
$\langle \overrightarrow{v}_9,\overrightarrow{v}_1\rangle=0, \langle \overrightarrow{v}_9,\overrightarrow{v}_2\rangle=250, \langle \overrightarrow{v}_9,\overrightarrow{v}_3\rangle=-200, \langle \overrightarrow{v}_9,\overrightarrow{v}_4\rangle=510, \langle \overrightarrow{v}_9,\overrightarrow{v}_5\rangle=0,\langle \overrightarrow{v}_9,\overrightarrow{v}_6\rangle=-460,\langle \overrightarrow{v}_9,\overrightarrow{v}_7\rangle=1960,\langle \overrightarrow{v}_9,\overrightarrow{v}_8\rangle=137300$.
\item {\bf V\'ertice $\mathbf v_{10}= (3,2)$}.\\
 Let $\overrightarrow{v}_{10} = k_{10,1}\overrightarrow{e}_1 + k_{10,2}\overrightarrow{e}_2 + k_{10,3}\overrightarrow{e}_3 + k_{10,4}\overrightarrow{e}_4 + k_{10,5}\overrightarrow{e}_5$.

El v\'ertice $v_{10}$ es adyacente  con  $v_1, v_3, v_4, v_5, v_7$, y $v_8$  en $\overline{K_3\square P_4}$ pero no es adyacente  con  $v_2,v_6$ ni  con  $v_9$.  De estas  condiciones obtenemos: $k_{10,1} + k_{10,4}-g_{10,1} = 0, k_{10,2}= 0, k_{10,1} + k_{10,3} - g_{10,3} = 0, k_{10,1} + k_{10,2} - k_{10,3} + k_{10,4} - g_{10,4} = 0,  -2 k_{10,1} + 2 k_{10,2} + k_{10,3} + 2 k_{10,4}-g_{10,5}= 0, 5 k_{10,1} + 2 k_{10,3} + 4 k_{10,4} = 0, 8 k_{10,1} + 4 k_{10,2} - 8 k_{10,3} - 6 k_{10,4} - g_{10,7} = 0, 13 k_{10,1} + 56 k_{10,2} + 53 k_{10,3} - 16 k_{10,4} - g_{10,8} = 0$, y $60k_{10,1}+250k_{10,2}-260k_{10,3}-60k_{10,4}=0$.
La matriz  del  sistema   es :
\newpage
{\tiny
$$
\left(
\begin{array}{ccccccccccc}
 1 & 0 & 0 & 1 & 0 & -1 & 0 & 0 & 0 & 0 & 0 \\
 0 & 1 & 0 & 0 & 0 & 0 & 0 & 0 & 0 & 0 & 0 \\
 1 & 0 & 1 & 0 & 0 & 0 & -1 & 0 & 0 & 0 & 0 \\
 1 & 1 & -1 & 1 & 0 & 0 & 0 & -1 & 0 & 0 & 0 \\
 -2 & 2 & 1 & 1 & 0 & 0 & 0 & 0 & -1 & 0 & 0 \\
 5 & 0 & 2 & 4 & 0 & 0 & 0 & 0 & 0 & 0 & 0 \\
 8 & 4 & -8 & -6 & 0 & 0 & 0 & 0 & 0 & -1 & 0 \\
 13 & 56 & 53 & -16 & 0 & 0 & 0 & 0 & 0 & 0 & -1 \\
 60 & 250 & -260 & -60 & 0 & 0 & 0 & 0 & 0 & 0 & 0
\end{array}
\right)\sim
\left(
\begin{array}{ccccccccccc}
 1 & 0 & 0 & 0 & 0 & 0 & 0 & 0 & 0 & 0 & -\frac{46}{3165} \\
 0 & 1 & 0 & 0 & 0 & 0 & 0 & 0 & 0 & 0 & 0 \\
 0 & 0 & 1 & 0 & 0 & 0 & 0 & 0 & 0 & 0 & -\frac{9}{1055} \\
 0 & 0 & 0 & 1 & 0 & 0 & 0 & 0 & 0 & 0 & \frac{71}{3165} \\
 0 & 0 & 0 & 0 & 0 & 1 & 0 & 0 & 0 & 0 & \frac{5}{633} \\
 0 & 0 & 0 & 0 & 0 & 0 & 1 & 0 & 0 & 0 & -\frac{73}{3165} \\
 0 & 0 & 0 & 0 & 0 & 0 & 0 & 1 & 0 & 0 & \frac{52}{3165} \\
 0 & 0 & 0 & 0 & 0 & 0 & 0 & 0 & 1 & 0 & \frac{136}{3165} \\
 0 & 0 & 0 & 0 & 0 & 0 & 0 & 0 & 0 & 1 & -\frac{578}{3165}
\end{array}
\right)
$$
}
Luego, si  $g_{10,8} = 3165, g_{10,7} = 578$ then $g_{10,5} = -136, g_{10,4} = -52, g_{10,3} = 73, g_{10,1} = -5$ todas las  condicines de adyacencia  y ortogonalidad se satisfacen. Tambi\'en  $k_{10,4} = -71, k_{10,3} = 27, k_{10,2} = 0$, and $k_{10,1} = 46$. Entonces ,
$$
\overrightarrow{v}_{10}= 46\overrightarrow{e}_1+27\overrightarrow{e}_3-71\overrightarrow{e}_4.
$$
Luego ,$
\langle \overrightarrow{v}_{10},\overrightarrow{v}_1\rangle=-25, \langle \overrightarrow{v}_{10},\overrightarrow{v}_2\rangle=0, \langle \overrightarrow{v}_{10},\overrightarrow{v}_3\rangle=73, \langle \overrightarrow{v}_{10},\overrightarrow{v}_4\rangle=-52, \langle \overrightarrow{v}_{10},\overrightarrow{v}_5\rangle=-207,\langle \overrightarrow{v}_{10},\overrightarrow{v}_6\rangle=0,\langle \overrightarrow{v}_{10},\overrightarrow{v}_7\rangle=578,\langle \overrightarrow{v}_{10},\overrightarrow{v}_8\rangle=3165,\langle \overrightarrow{v}_{10},\overrightarrow{v}_9\rangle=0$.
\item {\bf V\'ertice $\mathbf v_{11}= (3,3)$}.\\
Sea $\overrightarrow{v}_{11} = k_{11,1}\overrightarrow{e}_1 + k_{11,2}\overrightarrow{e}_2 + k_{11,3}\overrightarrow{e}_3 + k_{11,4}\overrightarrow{e}_4 + k_{11,5}\overrightarrow{e}_5$.

El v\'ertice $v_{11}$ es adyacente  con  $v_1, v_2, v_4, v_5, v_6,v_8$ y $v_9$ en  $\overline{K_3\square P_4}$ pero no es adyacente  con  $v_3, v_7$ y $v_{10}$.  De aqui  se o btienen las ecuaciones : $ k_{11,1} + k_{11,4} - g_{11,1} = 0, k_{11,2}-g_{11,2}=0, k_{11,1} + k_{11,3} = 0, k_{11,1} + k_{11,2} - k_{11,3} + k_{11,4} - g_{11,4} = 0, -2 k_{11,1} + 2 k_{11,2}+ k_{11,3} + 2 k_{11,4} - g_{11,5} = 0, 5 k_{11,1} + 2 k_{11,3} + 4 k_{11,4}-g_{11,6} = 0,  8 k_{11,1} + 4 k_{11,2} - 8 k_{11,3} - 6 k_{11,4} = 0, 13 k_{11,1} + 56 k_{11,2} + 53 k_{11,3} - 16 k_{11,4} - g_{11,8} =0, 60 k_{11,1} + 250 k_{11,2} - 260 k_{11,3} - 60 k_{11,4}-g_{11,9} = 0$, y $46k_{11,1}+27k_{11,3}-71k_{11,4}=0$. La matriz del  sistema  es:
{\tiny
$$
\left(
\begin{array}{cccccccccccc}
 1 & 0 & 0 & 1 & 0 & -1 & 0 & 0 & 0 & 0 & 0 & 0 \\
 0 & 1 & 0 & 0 & 0 & 0 & -1 & 0 & 0 & 0 & 0 & 0 \\
 1 & 0 & 1 & 0 & 0 & 0 & 0 & 0 & 0 & 0 & 0 & 0 \\
 1 & 1 & -1 & 1 & 0 & 0 & 0 & -1 & 0 & 0 & 0 & 0 \\
 -2 & 2 & 1 & 2 & 0 & 0 & 0 & 0 & -1 & 0 & 0 & 0 \\
 5 & 0 & 2 & 4 & 0 & 0 & 0 & 0 & 0 & -1 & 0 & 0 \\
 8 & 4 & -8 & -6 & 0 & 0 & 0 & 0 & 0 & 0 & 0 & 0 \\
 13 & 56 & 53 & -16 & 0 & 0 & 0 & 0 & 0 & 0 & -1 & 0 \\
 60 & 250 & -260 & -60 & 0 & 0 & 0 & 0 & 0 & 0 & 0 & -1 \\
 46 & 0 & 27 & -71 & 0 & 0 & 0 & 0 & 0 & 0 & 0 & 0
\end{array}
\right)\sim
\left(
\begin{array}{cccccccccccc}
 1 & 0 & 0 & 0 & 0 & 0 & 0 & 0 & 0 & 0 & 0 & \frac{71}{42295} \\
 0 & 1 & 0 & 0 & 0 & 0 & 0 & 0 & 0 & 0 & 0 & -\frac{511}{84590} \\
 0 & 0 & 1 & 0 & 0 & 0 & 0 & 0 & 0 & 0 & 0 & -\frac{71}{42295} \\
 0 & 0 & 0 & 1 & 0 & 0 & 0 & 0 & 0 & 0 & 0 & \frac{19}{42295} \\
 0 & 0 & 0 & 0 & 0 & 1 & 0 & 0 & 0 & 0 & 0 & \frac{18}{8459} \\
 0 & 0 & 0 & 0 & 0 & 0 & 1 & 0 & 0 & 0 & 0 & -\frac{511}{84590} \\
 0 & 0 & 0 & 0 & 0 & 0 & 0 & 1 & 0 & 0 & 0 & -\frac{189}{84590} \\
 0 & 0 & 0 & 0 & 0 & 0 & 0 & 0 & 1 & 0 & 0 & -\frac{686}{42295} \\
 0 & 0 & 0 & 0 & 0 & 0 & 0 & 0 & 0 & 1 & 0 & \frac{289}{42295} \\
 0 & 0 & 0 & 0 & 0 & 0 & 0 & 0 & 0 & 0 & 1 & -\frac{17452}{42295}
\end{array}
\right)
$$
}

Asi, si  $g_{11,9} = 84590$ then $ g_{11,8} = 34904, g_{11,6},g_{11,5} = 1372, g_{11,4} = 189, g_{11,2} = 511, g_{11,1} = -180$ aTodas las  condiciones de adyacencia   y  ortogonalidad  se satisfacen. Tambi\'en $k_{11,4} =-38, k_{11,3} =142, k_{11,2} = 511$, and $k_{11,1} = -142$. As\1,
$$
\overrightarrow{v}_{11}=-142\overrightarrow{e}_1+511\overrightarrow{e}_2+142\overrightarrow{e}_3-38\overrightarrow{e}_4.
$$
Luego,\
$\langle \overrightarrow{v}_{11},\overrightarrow{v}_1\rangle=-180, \langle \overrightarrow{v}_{11},\overrightarrow{v}_2\rangle=511, \langle \overrightarrow{v}_{11},\overrightarrow{v}_3\rangle=0, \langle \overrightarrow{v}_{11},\overrightarrow{v}_4\rangle=189,\\  \langle \overrightarrow{v}_{11},\overrightarrow{v}_5\rangle=1372,\langle \overrightarrow{v}_{11},\overrightarrow{v}_6\rangle=-578,\langle \overrightarrow{v}_{11},\overrightarrow{v}_7\rangle=0,\langle \overrightarrow{v}_{11},\overrightarrow{v}_8\rangle=34904,\langle \overrightarrow{v}_{11},\overrightarrow{v}_9\rangle=84590,\langle \overrightarrow{v}_{11}$, and $\overrightarrow{v}_{10}\rangle=0.
$
\item {\bf V\'ertice $\mathbf v_{12}= (3,4)$}.\\
  Sea $\overrightarrow{v}_{12} = k_{12,1}\overrightarrow{e}_1 + k_{12,2}\overrightarrow{e}_2 + k_{12,3}\overrightarrow{e}_3 + k_{12,4}\overrightarrow{e}_4 + k_{12,5}\overrightarrow{e}_5$.

El v\'ertice $v_{12}$ es adyacente  con  $v_1, v_2, v_3, v_5, v_6, v_7 , v_9$, y  $v_{10}$  en $\overline{K_3\square P_4}$  pero no es adyacente  con  $v_4, v_8$ and $v_{11}$.

de estas  condiciones se  obtienen las ecuaciones  $k_{12,1} + k_{12,4} - g_{12,1} = 0, k_{12,2} - g_{12,2} = 0, k_{12,1} + k_{12,3}-g_{12,3} = 0, k_{12,1} + k_{12,2} - k_{12,3} + k_{12,4} = 0, -2 k_{12,1} + 2 k_{12,2}+ k_{12,3} + 2 k_{12,4} - g_{12,5} = 0, 5 k_{12,1} + 2 k_{12,3} + 4 k_{12,4} - g_{12,6}= 0, 8 k_{12,1} + 4 k_{12,2} - 8 k_{12,3} - 6 k_{12,4}-g_{12,7} = 0, 13 k_{12,1} + 56 k_{12,2} + 53 k_{12,3} - 16 k_{11,4} = 0, 60 k_{12,1} + 250 k_{12,2} - 260 k_{11,3} - 60 k_{12,4} - g_{12,9} = 0, 46 k_{12,1} + 27 k_{12,3} - 71 k_{12,4}-g_{12,10}= 0$, and $-142 k_{12,1}+511 k_{12,2}+142k_{12,3}-38k_{12,4}=0$. La La  matriz  del sistema  es:
{\tiny
$$
\left(
\begin{array}{ccccccccccccc}
 1 & 0 & 0 & 1 & 0 & -1 & 0 & 0 & 0 & 0 & 0 & 0 & 0 \\
 0 & 1 & 0 & 0 & 0 & 0 & -1 & 0 & 0 & 0 & 0 & 0 & 0 \\
 1 & 0 & 1 & 0 & 0 & 0 & 0 & -1 & 0 & 0 & 0 & 0 & 0 \\
 1 & 1 & -1 & 1 & 0 & 0 & 0 & 0 & 0 & 0 & 0 & 0 & 0 \\
 -2 & 2 & 1 & 2 & 0 & 0 & 0 & 0 & -1 & 0 & 0 & 0 & 0 \\
 5 & 0 & 2 & 4 & 0 & 0 & 0 & 0 & 0 & -1 & 0 & 0 & 0 \\
 8 & 4 & -8 & -6 & 0 & 0 & 0 & 0 & 0 & 0 & -1 & 0 & 0 \\
 13 & 56 & 53 & -16 & 0 & 0 & 0 & 0 & 0 & 0 & 0 & 0 & 0 \\
 60 & 250 & -260 & -60 & 0 & 0 & 0 & 0 & 0 & 0 & 0 & -1 & 0 \\
 46 & 0 & 27 & -71 & 0 & 0 & 0 & 0 & 0 & 0 & 0 & 0 & -1 \\
 -142 & 511 & 142 & -38 & 0 & 0 & 0 & 0 & 0 & 0 & 0 & 0 & 0
\end{array}
\right)\sim
$$
$$
\left(
\begin{array}{ccccccccccccc}
 1 & 0 & 0 & 0 & 0 & 0 & 0 & 0 & 0 & 0 & 0 & 0 & -\frac{45}{10589} \\
 0 & 1 & 0 & 0 & 0 & 0 & 0 & 0 & 0 & 0 & 0 & 0 & -\frac{2288}{1005955} \\
 0 & 0 & 1 & 0 & 0 & 0 & 0 & 0 & 0 & 0 & 0 & 0 & \frac{7803}{1005955} \\
 0 & 0 & 0 & 1 & 0 & 0 & 0 & 0 & 0 & 0 & 0 & 0 & \frac{14366}{1005955} \\
 0 & 0 & 0 & 0 & 0 & 1 & 0 & 0 & 0 & 0 & 0 & 0 & \frac{10091}{1005955} \\
 0 & 0 & 0 & 0 & 0 & 0 & 1 & 0 & 0 & 0 & 0 & 0 & -\frac{2288}{1005955} \\
 0 & 0 & 0 & 0 & 0 & 0 & 0 & 1 & 0 & 0 & 0 & 0 & \frac{3528}{1005955} \\
 0 & 0 & 0 & 0 & 0 & 0 & 0 & 0 & 1 & 0 & 0 & 0 & \frac{40509}{1005955} \\
 0 & 0 & 0 & 0 & 0 & 0 & 0 & 0 & 0 & 1 & 0 & 0 & \frac{10339}{201191} \\
 0 & 0 & 0 & 0 & 0 & 0 & 0 & 0 & 0 & 0 & 1 & 0 & -\frac{191972}{1005955} \\
 0 & 0 & 0 & 0 & 0 & 0 & 0 & 0 & 0 & 0 & 0 & 1 & -\frac{743848}{201191}
\end{array}
\right)
$$
}

Asi, si  $g_{12,10} =1005955$ then $ g_{12,9} =3719240 , g_{12,7} =191972 , g_{12,6} =-51695 , g_{12,5}=-40509, g_{12,3} =-3528, g_{12,2} =2288, g_{12,1} =-10091$ se satisfacen   todas las  condiciones de adyacencia  y ortogonalidad. Tambi\'en  $k_{12,4} =-14366, k_{12,3} =-7803, k_{12,2} =2288$, y $k_{12,1} =4275$. En consecuencia,
$$
\overrightarrow{v}_{12}=4275 \overrightarrow{e}_1+2288 \overrightarrow{e}_2-7803 \overrightarrow{e}_3-14366\overrightarrow{e}_4.
$$
Luego,
$
\langle\overrightarrow{v}_{12},\overrightarrow{v}_1\rangle=-10091,\langle\overrightarrow{v}_{12},\overrightarrow{v}_2\rangle=2288, \langle\overrightarrow{v}_{12},\overrightarrow{v}_3\rangle=-3528, \langle\overrightarrow{v}_{12},\overrightarrow{v}_1\rangle=0, \langle\overrightarrow{v}_{12},\overrightarrow{v}_5\rangle=-40509,\langle\overrightarrow{v}_{12},\overrightarrow{v}_6\rangle=-51695,
\langle\overrightarrow{v}_{12},\overrightarrow{v}_7\rangle=191972,\langle\overrightarrow{v}_{12},\overrightarrow{v}_1\rangle=0,
\langle\overrightarrow{v}_{12},\overrightarrow{v}_9\rangle=3719240,\langle\overrightarrow{v}_{12},\overrightarrow{v}_{10}\rangle=1005955,
\langle\overrightarrow{v}_{12},\overrightarrow{v}_{11}\rangle=0
$
\end{itemize}
La  Matriz de Gram   $[G]_{i,j}=[\langle v_i,v_j\rangle]=0$  para la  representaci\'on  $\langle\overrightarrow{v}_{1},\cdots \langle\overrightarrow{v}_{12}$ es:
{\tiny
$$
\left(
\begin{array}{cccccccccccc}
 2 & 0 & 1 & 2 & 0 & 9 & 2 & -3 & 0 & -25 & -180 & -10091 \\
 0 & 1 & 0 & 1 & 2 & 0 & 4 & 56 & 250 & 0 & 511 & 2288 \\
 1 & 0 & 2 & 0 & -1 & 7 & 0 & 66 & -200 & 73 & 0 & -3528 \\
 2 & 1 & 0 & 4 & 1 & 7 & 14 & 0 & 510 & -52 & 189 & 0 \\
 0 & 2 & -1 & 1 & 13 & 0 & -28 & 107 & 0 & -207 & 1372 & -40509 \\
 9 & 0 & 7 & 7 & 0 & 45 & 0 & 107 & -460 & 0 & -578 & -51695 \\
 2 & 4 & 0 & 14 & -28 & 0 & 180 & 0 & 3920 & 578 & 0 & 191972 \\
 -3 & 56 & 66 & 0 & 107 & 107 & 0 & 6370 & 1960 & 3165 & 34904 & 0 \\
 0 & 250 & -200 & 510 & 0 & -460 & 3920 & 1960 & 137300 & 0 & 84590 & 3719240 \\
 -25 & 0 & 73 & -52 & -207 & 0 & 578 & 3165 & 0 & 7886 & 0 & 1005955 \\
 -180 & 511 & 0 & 189 & 1372 & -578 & 0 & 34904 & 84590 & 0 & 302893 & 0 \\
 -10091 & 2288 & -3528 & 0 & -40509 & -51695 & 191972 & 0 & 3719240 & 1005955 & 0 & 290779334
\end{array}
\right)
$$
}
Calculando los  valores propios de esta matriz usando MATHEMATICA 7,0  con  una aproximacion a  5  decimales  se obtienen los  valores $2.9083\times 10^8, 3.3586\times 10^5, 64461., 3167.1$ y $ 0$  lo  que indica son nonegativos    y por lo tanto el c\'alculo anterior  corresponde a  una   representaci\'on  ortogonal  de vectores linealmente independientes   del  grafo $K_3\square P_4$  y ademas la matriz  es semidefinida positiva.

En resumen  una representaci\'on  ortogonal   de $\overrightarrow{K_3 \square P_4}$ de vectores linealmente independientes  es :
\begin{itemize}
\item $\overrightarrow{v}_1 = \overrightarrow{e}_1 + \overrightarrow{e}_4$,
\item $\overrightarrow{v}_2 = \overrightarrow{e}_2$,
\item $\overrightarrow{v}_3 = \overrightarrow{e}_1 + \overrightarrow{e}_3$,
\item $\overrightarrow{v}_4 := \overrightarrow{e}_1 + \overrightarrow{e}_2 - \overrightarrow{e}_3 + \overrightarrow{e}_4$,
\item $\overrightarrow{v}_5 := -2 \overrightarrow{e}_1 + 2 \overrightarrow{e}_2 + \overrightarrow{e}_3 + 2 \overrightarrow{e}_4$,
\item $\overrightarrow{v}_6 := 5 \overrightarrow{e}_1 + 2 \overrightarrow{e}_3 + 4 \overrightarrow{e}_4$,
\item $\overrightarrow{v}_7 := 8 \overrightarrow{e}_1 + 4 \overrightarrow{e}_2 + -8 \overrightarrow{e}_3 - 6 \overrightarrow{e}_4$,
\item $\overrightarrow{v}_8 := 13 \overrightarrow{e}_1 + 56 \overrightarrow{e}_2 + 53 \overrightarrow{e}_3 - 16 \overrightarrow{e}_4$,
\item $\overrightarrow{v}_9 := 60 \overrightarrow{e}_1 + 250 \overrightarrow{e}_2 - 260\overrightarrow{ e}_3 - 60 \overrightarrow{e}_4$,
\item $\overrightarrow{v}_{10} := 46 \overrightarrow{e}_1 + 27 \overrightarrow{e}_3 - 71 \overrightarrow{e}_4$,
\item $\overrightarrow{v}_{11} := -142 \overrightarrow{e}_1 + 511 \overrightarrow{e}_2 + 142 \overrightarrow{e}_3 - 38\overrightarrow{ e}_4$,
\item $\overrightarrow{v}_{12} := 4275 \overrightarrow{e}_1 + 2288 \overrightarrow{e}_2 - 7803 \overrightarrow{e}_3 - 14366\overrightarrow{e}_4$.
\end{itemize}

\section{ Observaciones} 

El lector notar\'a   que en todas las matrices  usamos  $-1$ en las entradas de las  variables  $g_{m,j}, m=3,\dots, 12$ para obtener los  valores directamente.  En la demostracion  del  teorema  \ref{main} en \cite{PD} se usan variables  $(-g_{m,j}), m=3,\dots, 12$ para una mejor presentaci\'on del argumento.No obstante  el c\'alculo aqui hecho es equivalente.

Adem\'as  note    que   en  todas las matrices   las entradas de  la   quinta  columna  correspondiente a las coordenadas de   $\overrightarrow{e}_5$ are  all  zero.  esto significa   que  en  realidad  no necesitamos  calcular  la representaci\'on ortogonal en  $\Re^5$  sino que basta  calcularla  en  $\Re^4$. Sin  embargo,  no se puede establecer  a priori  si la  quinta  dimension  ser\1a  necesaria    y por lo tanto  es necesario calcular  vectores  en  $\Re^5$

Esta  t\'ecnica   aqui mostrada    se puede  aplicar  a  cualquier    grafo  simple  conexo  con lo  cual se tiene una herramenta  para  calcular  una  cota  superior  del  rango m\1nimo  semidefinido  de estos   grafos   y posiblemente   se podr\1a   determinar  si  la  conjetura delta  se satisface  para  un  grafo  simple  conexo  particular  es v\'alida. 

\section{Conclusi\'on}

La t\'ecnica  de encrustamientos  circulares para la  representaci\'on de  grafos  simples  y conexos  es una  t\'ecnica  que  se puede  aplicar   a  cualquier grafo  simple  conexo.  Es particularmente \'util   para  calcular las  representaciones ortogonales de   grafos  simples conexos  que a su  vez es esencial   en el estudio  de  problemas abiertos  en  an\'alisis matricial  como la  conjetura delta  y la  conjetura del   grafo  complementario. Por  la  forma  en que  se simplifica  la observaci\'on de  un  grafo simple conexo   y  sus   relaciones de  simetr\1a  encontradas  una  vez  que se  ha  aplicado la t\'ecnica de encrustamiento se tiene  una herramienta muy  util para  el estudio  de los  grafos  simples  conexos.

\section{Agradecimientos }
Deseo  agradecer a mi  director de tesis doctoral   Dr. Sivaram Narayan  por sus consejos  y  sugerencias    durante la  investigaci\'on    doctoral   sobre    $\delta$-grafos    por la cual  fue posible la escritura de este art\1culo.  Tambien  quiero agradecer a la Universidad de Costa  Rica y  particularmente  al  Centro de  Investigaciones  Matem\'aticas y  Metamatem\'aticas (CIMM)  por  el soporte econ\'omico para la investigaci\'on. Finalmente deseo  agradecer al departamento  de  matem\'atica de  Central Michigan University   por el  soporte econ\'omico  otorgado durante  mi  formaci\'on  doctoral.
\input{bibliography}
\end{document}

%% file: bibliography.tex
\renewcommand{\baselinestretch}{1} 

%% file: Encrustamientoscirculares.bbl
\small\normalsize
\begin{thebibliography}{99}
\bibitem{FW2}
AIM Minimum Rank-Special Graphs Work Group (Francesco Barioli, Wayne Barrett, Steven Butler, Sebastian M. Cioaba, Shaun M. Fallat, Chris Godsil, Willem Haemers, Leslie Hogben, Rana Mikkelson, Sivaram Narayan, Olga Pryporova, Irene Sciriha, Dragan Stevanovic, Hein Van Der Holst, Kevin Van Der Meulen, and Amy Wangsness), {\it Zero Forcing Sets and the Minimum Rank of Graphs}, Linear Algebra and its Applications, 428 (2008) 1628-1648.
\bibitem{FW}
Francesco Barioli, Wayne Barrett, Shaun M. Fallat, H. Tracy Hall, Leslie Hogben,  and Hein van der Holst, {\it On the Graph Complement Conjecture for Minimum Rank}, Linear Algebra and its Applications, in press, doi:10.1016/j.laa.2010.12.024.
\bibitem{FW1}
Francesco Barioli, Wayne Barrett, Shaun M. Fallat, H. Tracy Hall, Leslie Hogben, Bryan Shader, P. van den Driessche, and Hein van der Holst, {\it Zero Forcing Parameters and Minimum Rank Problems}, Linear Algebra and its Applications, 433 (2010) 401-411.
\bibitem{FS}
Francesco Barioli, Shaun M. Fallat, Lon H. Mitchell, and Sivaram Narayan, {\it Minimum Semidefinite Rank of Outerplanar Graphs and the Tree Cover Number}, Electronic Journal of Linear Algebra, 22 (2011) 10-21.
\bibitem{WH}
Wayne Barrett, Hein Van Der Holst, and Raphael Loewy, {\it Graphs Whose Minimal Rank is Two}, Electronic Journal of Linear Algebra, 11 (2004) 258-280.
\bibitem{BN}
Jonathan Beagley, Sivaram Narayan, Eileen Radzwion, Sara Rimer, Rachel  Tomasino, Jennifer Wolfe, and Andrew Zimmer , {\it On the Minimum Semidefinite  Rank  of a Graph  Using  Vertex Sums, Graphs  with $\msr(G)=|G|-2$, and  the $msr$s of  Certain  Graphs Classes },in NSF-REU Report  from Central Michigan University (Summer 2007).
\bibitem{AB}
Avi Berman, Shmuel Friedland, Leslie Hogben, Uriel G. Rothblum, and Bryan Shader, {\it An Upper Bound for the Minimum Rank of a Graph}, Linear Algebra and its Applications, 429 (2008) 1629-1638.
\bibitem{BO}
B\'ela Bollob\'as. \emph{Modern Graph Theory},Springer.Memphis, TN, 1998.
\bibitem{BM}
John Adrian Bondy, Uppaluri Siva Ramachandra  Murty, \emph{ Graph Theory},Springer.San  Francisco, CA, 2008.
\bibitem{MP2}
 Matthew Booth, Philip Hackney, Benjamin Harris, Charles Johnson, Margaret Lay, Terry Lenker, Lon H. Mitchell, Sivaram K. Narayan, Amanda Pascoe, and Brian D. Sutton, {\it On the Minimum Semidefinite Rank of a Simple Graph}, Linear and Multilinear Algebra, 59 (2011) 483-506.
\bibitem{MP3}
Matthew Booth, Philip Hackney, Benjamin Harris, Charles R. Johnson, Margaret Lay, Lon H. Mitchell, Sivaram K. Narayan, Amanda Pascoe, Kelly Steinmetz, Brian D. Sutton, and Wendy Wang, {\it On the Minimum Rank Among Positive Semidefinite Matrices with a Given Graph}, SIAM Journal on Matrix Analysis and Applications, 30 (2008) 731-740.
\bibitem{BLS}
 Andreas Brandstädt, Van Bang  Le, Jeremy Spinrad,{\it Graph Classes: A Survey, SIAM Monographs on Discrete Mathematics and Applications}, ISBN 0-89871-432-X.(1999) p 169.
\bibitem{RB1}
Richard Brualdi, Leslie Hogben, Bryan Shader, {\it AIM Workshop Spectra of Families of Matrices described by  Graphs, Digraphs,  and  Sign Patterns Final  Report: Mathematical  Results (Revised*)} Available at http://www.aimath.org/pastworkshop/matrixspectrumrep.pdf, (2007).
\bibitem{CH}
Gary Chartrand, Linda Lesniak, and Ping Zhang. \emph{ Graphs \& Digraphs},Taylor \&  Francis Group. Boca Raton, FL, 2011.
\bibitem{PD}Pedro,Diaz. {\it On the delta Conjecture  and  the Graph Complement Conjecture for Minimum Semidefinite rank of a  Graph}, Ph.D Dissertation, Central Michigan University, July (2014).
\bibitem{JE}
Jason Ekstrand, Craig Erickson, H. Tracy Hall, Diana Hay, Leslie Hogben, Ryan Johnson, Nicole Kingsley, Steven Osborne, Travis Peters, Jolie Roat, Arianne Ross, Darren D. Row, Nathan Warnberg, and Michael Young, {\it Positive Semidefinite Zero Forcing}, Linear Algebra and its Applications, 439 (2013) 1862-1874.
\bibitem{SL}
Shaun M. Fallat,  Leslie Hogben, {\it The Minimum Rank of Symmetric Matrices Described by a Graph: A Survey}, Linear Algebra and its Applications, 426 (2007) 558-582.
\bibitem{PB}
Philip Hackney, Benjamin Harris, Margaret Lay, Lon H. Mitchell, Sivaram K. Narayan, and Amanda Pascoe, {\it Linearly Independent Vertices and Minimum Semidefinite Rank}, Linear Algebra and its Applications, 431 (2009) 1105-1115.
\bibitem{VH}
Hein van der Holst, {\it Graphs whose Positive Semidefinite Matrices have Nullity at Most Two}, Linear Algebra and its Applications, 375 (2003) 1-11.
\bibitem{RC}
Roger Horn, Charles  Johnson, {\it Matrix Analysis}, Cambridge University Press, 1985.
\bibitem{LH1}
Leslie Hogben, {\it Minimum Rank Problems}, Linear Algebra and its Applications, 432 (2009) 1961-1974.
\bibitem{LH2}
Leslie Hogben, {\it Orthogonal Representations, Minimum Rank, and Graph Complements}, Linear Algebra and its Applications, 428 (2008) 2560-2568.
\bibitem{YL}
Yunjiang Jiang, Lon H. Mitchell, and Sivaram K. Narayan, {\it Unitary Matrix Digraphs and Minimum Semidefinite Rank}, Linear Algebra and its Applications, 428 (2008) 1685-1695.
\bibitem{LM}
Lon  Mitchell, {\it On the Graph Complement Conjecture for Minimum Semidefinite Rank}, Linear Algebra and its Applications (2011), in press, doi:10.1016/j.laa.2011.03.011.
\bibitem{LSA2}
Lon Mitchell, Sivaram K. Narayan, and Andrew M. Zimmer, {\it Lower Bounds in Minimum Rank Problems}, Linear Algebra and its Applications, 432 (2010) 430-440.
\bibitem{SY1}
Sivaram K. Narayan, Yousra Sharawi, {\it Bounds on Mi\-ni\-mum Se\-mi\-de\-fi\-ni\-te Rank of Graphs}, Linear and Multilinear Algebra (2014). Retrieved from\\
http://dx.doi.org/10.1080/03081087.2014.898763 on  June 13, 2014.
\bibitem{PN}
Peter  Nylen, {\it Minimum-Rank Matrices with Prescribed Graph}, Linear Algebra and its Applications, 248 (1996) 303-316.
\bibitem{TP}
Travis Peters, {\it Positive Semidefinite Maximum Nullity and Zero Forcing Number}, Electronic Journal of Linear Algebra,23 (2012)828-829.
\bibitem{RW}
Ronald Read, Robin Wilson.\ \emph{An  Atlas of  Graphs}, Oxford University Press., 1998, 125-305.
\bibitem{SY}
Yousra Sharawi, {\it  Minimum Semidefinite  Rank of  a Graph}.Ph.D. Dissertation. Central  MIchigan  University, December (2011).
\bibitem{DW}
Douglas B. West, \emph{Introduction to Graph Theory}, Prentice Hall Inc., 1996.
\end{thebibliography}
